\documentclass[11pt]{amsart}

\usepackage{amsbsy,amssymb,amscd,amsfonts,latexsym,amstext,delarray,
amsmath,epsfig,graphicx,bm, url}
\usepackage[all]{xy}
 \usepackage[mathscr]{eucal} 
 \usepackage{color}
\input xypic
\usepackage{hyperref}

\newtheorem{thm}{Theorem}[section]
\newtheorem{prop}[thm]{Proposition}
\newtheorem{cor}[thm]{Corollary}
\newtheorem{lem}[thm]{Lemma}
\newtheorem{defn}[thm]{Definition}
\newtheorem{rem}[thm]{Remark}
\newtheorem{ex}[thm]{Example}

\numberwithin{equation}{section}

\def\C{{\mathbb C}}

\def\N{{\mathbb N}}

\def\Z{{\mathbb Z}}
\def\R{{\mathbb R}}

\def\fH{{\mathfrak H}}

\def\cutint{{\int_D \!\!\!\!\!\!\!-}}

\def\cA{{\mathcal A}}
\def\cB{{\mathcal B}}

\newcommand{\sC}{{\mathscr C}}
\def\cD{{\mathcal D}}

\def\cF{{\mathcal F}}
\def\cG{{\mathcal G}}

\def\cK{{\mathcal K}}
\def\cL{{\mathcal L}}

\def\cP{{\mathcal P}}

\def\cR{{\mathcal R}}
\def\cS{{\mathcal S}}
\newcommand{\sS}{{\mathscr S}}

\def\cU{{\mathcal U}}

\newcommand{\ie}{{\it i.e.\/}\ }
\newcommand{\eg}{{\it e.g.\/}\ }
\newcommand{\cf}{{\it cf.\/}\ }
\newcommand{\comp}{{\it comp.\/}\ }
\newcommand{\opcit}{{\it op.cit.\/}\ }

\def\text{\hbox}

\def\Aut{{\rm Aut}}
\def\Isom{{\rm Isom}}
\def\Sim{{\rm Sim}}

\def\Dom{{\rm Dom}}

\def\End{{\rm End}}

\def\id{{\rm Id}}
\def\prs{\s_{\rm pr}}

\def\Ind{{\rm Ind}}
\def\Index{{\rm Index}}
\def\Ker{{\rm Ker}}

\def\Res{{\rm Res}}

\def\Tr{{\rm Tr}}

\def\vol{{\rm vol}}

\def\a{\alpha}
\def\b{\beta}

\def\d{\delta}
\def\e{\epsilon}
\def\g{\gamma}

\def\i{\iota}
\def\lb{\lambda}
\def\om{\omega}
\def\s{\sigma}

\def\ve{\varepsilon}

\def\z{\zeta}

\def\D{\Delta}
\def\G{\Gamma}

\def\nb{\Nabla}
\def\Om{\Omega}

\def\fl{\forall}

\def\nb{\nabla}

\def\ot{\otimes}
\def\hot{\hat\otimes}
\def\part{\partial}

\def\sbs{\subset}

\def\ra{\rightarrow}

\def\text{\hbox}

\def\bash{\backslash}

\def\fl{\forall}

\def\nb{\nabla}
\def\ot{\otimes}

\def\ra{\rightarrow}

\def\sbs{\subset}

\def\trr{\triangleright}

\def\wt{\widetilde}
\def\td{\tilde}

\newcommand{\bigsetdef}[2]{\bigl\{ #1 \,\bigm|\, #2\bigr\}}

\newcommand{\qqq}{{\,,\quad \forall\,}}

\def\Dirac{{D \hspace{-7pt} \slash}}
\def\Jb{{J \hspace{-6pt} \bash}}

\begin{document}

\title{Local index formula and twisted spectral triples}
 \author{Henri Moscovici}
\address{Department of mathematics \\
The Ohio State University \\
 Columbus, OH 43210, USA}
\email{henri@math.ohio-state.edu}

\thanks{This work was
 supported by the National Science Foundation
awards nos. DMS-0245481  and DMS-0652167}

\begin{abstract}  
We prove a local index formula for a class of twisted spectral
triples of type III modeled on the transverse geometry
of conformal foliations with locally constant 
 transverse conformal factor. Compared with the earlier proofs in
 the untwisted case, the novel aspect resides in the fact that
 the twisted analogues of the JLO entire cocycle and of its retraction
are no longer cocycles in their respective Connes bicomplexes. 
We show however that the passage to the
infinite temperature limit, respectively the integration along the full temperature range
against the Haar measure of the positive half-line,
has the remarkable effect of curing in both cases the deviations from the cocycle 
and transgression identities. 
\end{abstract}

\maketitle

\centerline{\it Dedicated to Alain Connes, with admiration, friendship and gratitude}

\medskip

\section*{Introduction}

The local-global principle, as epitomized by the Atiyah-Singer index theorem
but in the larger operator theoretic framework,
has played a pivotal role in Alain Connes' overarching design of the foundations of
Noncommutative Geometry. Although I was too overwhelmed by his brilliant intellect 
and fantastic mathematical insight to fully realize it at the time,
this very theme was in fact the subtext of
our first mathematical conversation, in the autumn of 1978, while we were both visiting
the Institute for Advanced Study. As his program advanced, the theme
gradually evolved into a perennial context for a
substantial part of our collaborative work and, last but not least, it became a
pretext for a close, lifelong friendship. As a token of my deep appreciation,
I thought it would befit the occasion to try to run anew the machinery that has emerged,
this time in the presence of a twist.
  
\medskip

The basic template for a space in the Connes program is encoded in
the notion of spectral triple. In our recent joint paper~\cite{C-M4} we showed
that this notion can be adapted to include certain type III spaces 
by the simple device 
of incorporating in it a twisting 
automorphism of the algebra of coordinates. 
The paradigmatic examples of such spaces are those
describing the transverse geometry of a foliation of  codimension $1$,
or of a conformal foliation of arbitrary codimension. While it is always possible to
associate a spectral triple of type II to any foliation by
passing to the frame bundle of a complete transversal (\cf~\cite{C-M2}),
the very construction introduces a large number of additional parameters, which
in turn makes
the task of computing the characteristic classes quite formidable
(see~\cite{C-M3}). The twisting device on the other hand, whenever applicable, allows 
to bypass the extra step of geometric reduction to type II.
\medskip

Since the primary effect of the twisting automorphism
is the replacement of the bimodule of
noncommutative differential forms of a spectral triple
with a bimodule of twisted differential forms, one would have
expected the characteristic classes of a twisted spectral triple to be captured by a twisted
version of the Connes-Chern character, landing in twisted cyclic cohomology.
Somewhat counterintuitively, it turned out that no cohomological twisting is
actually needed, and that Connes' original construction of the 
Chern character in noncommutative geometry~\cite{ncdg}  remains in fact 
operative in the twisted case as well.  
The natural question that arises is whether the Connes-Chern character
of a twisted spectral triple can also be expressed in local terms as in~\cite{C-M2}, 
by means of a residue integral that eliminates all quantum infinitesimal
perturbations of order strictly larger than $1$. 
\medskip

In this paper we produce such a local index formula for a class of spectral
triples twisted by scaling automorphisms, modeled on the geometry of a conformal foliation 
whose holonomy consists of germs of conformal transformations of $\R^n$. 
The proof is patterned on the
strategy that evolved over a number of years of joint (and joyful)
work leading up to the residue index formula~\cite{C-M0, C-M, C-M1, C-M2},
 with the important difference that the twisted counterpart of the
JLO cocycle~\cite{JLO}, which played a key role in our earlier proof (as well as in
in Higson's~\cite{Higson}) is no longer a cocycle.  
\medskip

The plan of the paper is as follows. In \S 1 we recall from~\cite{C-M4}
 the basic definitions concerning twisted spectral triples and their characters.
 Extrapolating from the expression of local Hochschild cocycle in \opcit we then
 make a straightforward Ansatz in \S 2, predicting the form of the twisted local 
 formula for the Connes-Chern
 character. In \S 3 we test the Ansatz on ``real-life''
 examples of twisted spectral triples occurring in conformal geometry. These are
 the spectral triples describing the transverse geometry of the conformal foliations
 whose holonomy is given by germs of M\"obius transformations of $S^n$. We
 conclude that the Ansatz holds true if the holonomy is restricted to
 the parabolic subgroup preserving a point, or equivalently to the similarity 
 transformations of $R^n$. 
 \medskip

 The main results of the paper are proved in \S4,
 where we establish the validity of the Ansatz for an abstract class of twisted spectral triples,
 modeled on the conformal foliations with locally constant 
 transverse conformal factor. While  
the twisted entire cochain analogous to \cite{JLO}
is no longer a cocycle, passing to the
infinite temperature limit has the remarkable effect of
restoring the cocycle identity, and the resulting cocycle can
be expressed in terms of a residue integral as in~\cite{C-M2}. To show that
this residue cocycle represents the Connes-Chern character, one needs to
transit through a transgressed version, as in~\cite{C-M1}. In turn,
the transgression process does yield a genuine
cocycle because it involves integrating along the full temperature range, 
with respect to the Haar measure of $\R^+$,
which miraculously cures again the deviation from the cocycle identity.

 \bigskip

\section{Twisted spectral triples and their characters} 

We begin by briefly reviewing the notion of
twisted spectral triple of type III, introduced
in~\cite{C-M4}, together with some of its basic properties. 
 
\subsection{Twisted spectral triple} A {\it twisted spectral
triple} $(\cA , \fH , D, \s)$ consists of a local Banach $\ast$-algebra 
 $\cA$ represented in the Hilbert space $\fH$ by bounded operators, 
an automorphism  $\sigma \in \Aut (\cA)$, and 
 a self-adjoint
unbounded operator $D$ such that, for any $a \in \cA$, 
\begin{eqnarray} \label{compop}
&& a (D^2 + 1)^{-\frac{1}{2}} \in \cK (\fH) \quad \text{(compact operators)}; \\ \label{boundedcom}
&& a (\Dom D)\, \subset \, \Dom D \quad \text{and} \quad
[D, a]_\s:=D\,a\, -\,\sigma(a)\,D \quad \text{is bounded} ; \\ \label{sigmainv}
&& \sigma(a^*)=\,(\sigma^{-1}(a))^* .
\end{eqnarray}
 
A $\Z_2$-{\it graded} (or {\it even}) $\sigma$-spectral triple has the additional
datum of a grading operator
$$ \gamma = \gamma^* \in \cL (\fH) \, , \quad \gamma^2 = I
$$
which anticommutes with $D$, and commutes with the action of
$\cA$. In case the algebra itself is $\Z_2$-graded, the commutators
properties (including the twisted commutators) are understood in the
graded sense.

We shall be concerned with finitely-summable twisted spectral triples, \ie with
 those
that satisfy a stronger version of \eqref{compop}, namely the 
{\it $(p, \infty)$-summability} condition
\begin{equation} \label{nsumm}
 a (D^2 + 1)^{-\frac{1}{2}} \in \cL^{(p, \infty)} (\fH) , \qquad \forall \, a \in \cA ,
 \end{equation}
for some $\, 1 \leq p < \infty$ of the same parity as the spectral triple. The notation is
that of \cite[IV, 2.$\alpha$]{book} ,
\begin{eqnarray} \label{Lpinfty}
\cL^{(p, \infty)} (\fH) &=& \{ T \in \cK (\fH) \, ;  \quad \sum_{i=0}^N \mu_i (T) = O (N^{1-\frac{1}{p}}) \} ,
\qquad \text{if} \quad p > 1 , \\ \label{Lp1}
 \cL^{(1, \infty)} (\fH) &=& \{ T \in \cK (\fH) \, ;  \quad \sum_{i=0}^N \mu_i (T) = O (\log N) \} .
\end{eqnarray}

\subsection{Graded double} \label{grdbl}
As in the untwisted situation, there is a canonical way (\cf~\cite[Part I, \S 7]{ncdg})
to pass from an {\it ungraded} (or {\it odd}) twisted  spectral triple $(\cA , \fH , D, \s)$
to a $\Z_2$-graded twisted one over
the $\Z_2$-graded algebra
$$
\cA_{\rm gr} \, = \, \cA \ot C_1 .
$$ 
Here $C_1$ denotes the Clifford algebra 
${\rm Cliff} (\R) \ot \C$; its even part is $C_1^+ = \C \,1$, 
with $1$ the unit of $C_1$, while the odd part is $C_1^- = \C \, \e$, with $\e^2 =1$.
The automorphism remains $\s$, identified to $\s \ot \id$.
One constructs an $\cA_{\rm gr}$-module by first
letting $\fH_1 = \fH^+_1 \oplus  \fH^-_1$ be the
$\Z_2$-graded Hilbert space with $ \fH^{\pm}_1 = \C$ on which $C_1$ acts
via
$$
\lb 1 + \mu \e \mapsto  \begin{pmatrix} \lb& &\mu \cr
 \mu& &\lb 
  \end{pmatrix} , \qquad \lb, \mu \in \C ,
$$
and then taking 
$$
\fH_{\rm gr} = \fH \ot \fH_1, \quad \text{with} \quad  \fH_{\rm gr}^\pm = \fH \ot \fH_1^\pm ,
$$
on which $\cA_{\rm gr}$ acts via the exterior tensor product representation;
the corresponding grading operator is
$$ \g \, = \, \id_\fH \ot \g_1 , \quad \text{where} \quad
\g_1 = \begin{pmatrix} 1& &0 \cr
 0& &-1 
  \end{pmatrix} .
$$ 
Finally, as operator one takes
$$
D_{\rm gr}  \, = \, D \ot P_1 , \quad \text{where} \quad
P_1 = \begin{pmatrix} \, 0& &i \cr
 -i& &0 \end{pmatrix} :  \fH_1 \ra \fH_1 .
$$

\subsection{Invertible double} \label{invdbl}
When $D$ is not invertible, we shall resort to the construction 
described in~\cite[Part I, \S 6]{ncdg} (akin to the passage from the
Dirac operator on  flat space to the Dirac Hamiltonian with mass)
 to canonically associate to  $(\cA , \fH , D, \g , \s)$ a new
$\sigma$-spectral triple  $(\cA , \wt\fH , \wt{D}, \wt\g , \s)$ with invertible operator,
defined as follows.
With $\fH_1$ as above, one takes 
\begin{eqnarray*}
\wt\fH &=& \fH \hot \fH_1 \quad \text{(graded tensor product)} , \qquad
\wt\g = \g \hot \g_1 , \\
\wt{D} &=& D \hot \id + \id \hot F_1, \quad \text{where} \quad
F_1 = \begin{pmatrix} \, 0& &1 \cr
 1& &0 \end{pmatrix} :  \fH_1 \ra \fH_1 .
\end{eqnarray*}
The algebra $\cA$ is made to act on $\wt\fH$ via the representation
$$ a \in \cA \mapsto \tilde{a}:= a \hot e_{1} , \quad \text{where} \quad
e_{1} = \begin{pmatrix} \, 1& &0 \cr
 0& &0 \end{pmatrix} :  \fH_1 \ra \fH_1 .
$$

\subsection{Lipschitz regularity}   \label{untwist}
In \cite{C-M4} a
twisted spectral triple $(\cA , \fH , D, \s)$
 with invertible $D$ was called {\it Lipschitz regular} if it
satisfies the additional condition
\begin{equation}
[ |D|, a]_\s := |D|\,a-\,\sigma(a)\,|D| \quad  \hbox{is bounded}, \quad   \forall \, a \in \cA .
\label{boundedcom5}
\end{equation}
Such a twisted spectral triple can be `untwisted'  by passage to its
`phase' operator $\, F =  D \, |D|^{-1} $. Indeed, 
\begin{equation*} 
[F, a] \, = \, |D|^{-1} \, \big(  (D\,a\, -\,\sigma(a)\,D) \, -
\, (|D|\,a-\,\sigma(a)\,|D| )\, F \big)  ,
\end{equation*}
which shows that these commutators are compact operators,  
of the same order of magnitude as $\, D^{-1}$. Thus,  $\, (\fH , F)$
is a Fredholm module over $\cA$, defining
a $K_\ast$-cycle over the norm closure 
$C^*$-algebra of $\cA$. Moreover, if $(\cA , \fH , D, \s)$ is
$(p, \infty)$-summable so is $( \fH , F)$.

 We extend the Lipschitz regularity condition to a general 
 spectral triple $(\cA , \fH , D, \s)$
 by requiring that its invertible double $(\cA , \wt\fH , \wt{D}, \wt\g , \s)$ be 
 Lipschitz regular. Since
 \begin{equation*} 
 \wt{F} := \wt{D} \vert  \wt{D} \vert^{-1} = D (D^2 + 1)^{-\frac{1}{2}} \hot \id + 
 (D^2 + 1)^{-\frac{1}{2}} \hot F_1 ,
 \end{equation*}
 and therefore
  \begin{eqnarray*} 
 [\wt{F}, a]  \, =  \,
 [D (D^2 + 1)^{-\frac{1}{2}} , a] \hot e_{11} \, + \, 
  (D^2 + 1)^{-\frac{1}{2}} a\hot  e_{21} \, - \, a (D^2 + 1)^{-\frac{1}{2}} \hot e_{12} ,
 \end{eqnarray*}
Lipschitz regularity can be alternatively phrased as the requirement
\begin{equation}  \label{Lreg}
 [D (D^2 + 1)^{-\frac{1}{2}} , a]  \in \cK (\fH) , \qquad \fl \, a \in \cA \, .
\end{equation}
In the $(p, \infty)$-summable case these commutators belong to the ideal
$\cL^{(p, \infty)} (\fH)$ of $ \cK (\fH)$. 

\medskip

\subsection{Connes-Chern character} \label{chernchar}
By resorting if necessary
to the doubling procedures described in \ref{grdbl} and \ref{invdbl}, 
we may assume without essential loss of generality that the
$(p, \infty)$-summable twisted spectral triple $(\cA , \fH , D, \s)$
under consideration
is $\Z_2$-graded and $D$ invertible. 
We shall often do so in the sequel without any further mention.

Let $(\cA , \fH , D, \s)$ be such a twisted spectral triple which
is also {\it Lipschitz-regular}. Then,
as remarked in \ref{untwist}, it gives rise to a `phase' Fredholm module $( \fH , F)$ 
over $\, \cA$. In turn, the latter has a
well-defined Connes-Chern character in cyclic cohomology, \cf~\cite{ncdg}, 
which up to a normalizing constant is represented by the cyclic cocycle
\begin{equation}\label{CCh}
\tau^p_F (a_0 , a_1, \ldots , a_p) := 
\Tr \, (\gamma\, F\,  [F, a_0] \, [F, a_1] \cdots \, [F, a_p]) \, ,
\qquad  a_i \in \cA .
\end{equation}

In \cite{C-M4}, we have actually shown that a $(p, \infty)$-summable
twisted spectral triple as above admits a
Connes-Chern character without assuming Lipschitz regularity. Indeed, we showed
that the $(p+1)$-linear form on $\cA$
\begin{equation}\label{Cycocycle}
\tau^p_{F_\sigma} (a_0 , a_1, \ldots , a_p)\, :=  \,
\Tr \, (\gamma\, D^{-1}  [D, a_0]_\s \, D^{-1} [D, a_1]_\s \cdots \, D^{-1} [D, a_p]_\s),
\end{equation}
is a cyclic cocycle in $Z_\lambda^n (\cA)$, by means of which one can
recover the index pairing of $D$ with $K^\ast (\cA)$. More precisely, up to a
universal constant factor $c_p$, $\Index (\s (e) \, D \, e)$ is given by 
$\tau^p_{F_\sigma} (e , \ldots , e)$, for any class $[e] \in K^0 (\cA)$.

\subsection{Conformally perturbed spectral triple} \label{confpert}
An instructive class of
examples of  twisted spectral triples arises from
 conformal-type perturbations of ordinary spectral triples. 
Let  $(\cA , \fH , D) $ be a $(p, \infty)$-summable spectral triple,
and let $h = h^* \in \cA$. By setting
\begin{equation*}
D_h \,=\,e^{h}\,D\,e^{h} , \quad \text{and} \quad
 \sigma_h(a)=\,e^{2h}\,a\,e^{-2h}, \quad  \forall \, a \in \cA 
\end{equation*}
 one easily sees that
\begin{equation} \label{bound'}
 D_h\,a\, -\,\sigma_h(a)\,D_h  \, = \, e^{h} \, [D, \, \s_{h/2} (a)] \, e^{h}
\in \cL(\fH) \,,
\end{equation}
thus giving rise to a twisted spectral triple $(\cA , \fH , D_h, \s_h) $.
Noting that
\begin{equation}  \label{clue}
D_h^{-1} \, [D_h, a]_\s \, = \, e^{-h} \, D^{-1} \, [D, \, \s_{h/2} (a)] \, e^h \, ,
\end{equation}
one obtains the identity
\begin{eqnarray} \label{Ch}
&&\Tr \, (\g\,D_h^{-1} \, [D_h, a_0]_\s  \, D_h^{-1} \, [D_h, a_1]_\s \cdots \, 
D_h^{-1} \, [D_h, a_p]_\s ) \\ \notag
&&\qquad \qquad = \,  \Tr \, (\gamma\, D^{-1} [D,  \s_{h/2} (a_0)] \, D^{-1} [D,  \s_{h/2} (a_1)]
 \cdots \,D^{-1} [D,  \s_{h/2} (a_p)]) \, .
\end{eqnarray}
The right hand side is a cyclic cocycle on $\cA$ that
represents, for $h=1$ and up to normalization, the Connes-Chern character 
$$Ch^* (\cA , \fH , D) \in HC^* (\cA)
$$ 
of $(\cA , \fH , D)$ viewed as a Fredholm module, \cf~\cite[Part I, \S 6]{ncdg}.
It follows that the left hand side, which is
the cyclic cocycle obtained by composition with the inner automorphism $ \s_{h/2}$,
determines the same class in periodic
cyclic cohomology. This justifies regarding \eqref{Ch} as defining
the {\it periodic Connes-Chern character}
of the conformally perturbed spectral triple:
\begin{equation} \label{chpert}
 Ch^* (\cA , \fH , D_h, \s_h) \, := \,  Ch^* (\cA , \fH , D) \, \in \, HP^* (\cA) .
\end{equation}

\subsection{Local Hochschild cocycle} \label{hoch}
With the goal of extending
the local index formula of~\cite{C-M2} to twisted spectral triple, we 
looked in~\cite{C-M4} for 
an analogue of Connes' residue formula~\cite[IV.2.$\g$]{book} for
the Hochschild class of the Connes-Chern character. 

We recall that if  $(\cA , \fH , D)$ is an (even, invertible) 
$(p, \infty)$-summable spectral triple satisfying
the {\it smoothness condition}
\begin{equation} \label{R}
 \cA , \,  [D, \cA] \subset \bigcap_{k > 0} \Dom (\d^k) , \quad
\text{where} \quad \d (T) := [|D|, T ] \, ,
\end{equation}
the Hochschild cohomology class
$I \big(Ch^* (\cA , \fH , D)\big) \in HH^* (\cA)$ 
admits a {\it local} representation, 
given by the formula
\begin{equation}\label{Hcocycle1}
\varkappa_D (a_0 , a_1, \ldots , a_p) :=
\Tr_\om \big(\gamma a_0  [D, a_1] \cdots  [D, a_p] \, D^{-p}\big) ,
\qquad  a_i  \in \cA \, ,
\end{equation}
which defines a Hochschild cocycle.
Here $ \,\Tr_\om$ stands for a Dixmier trace (see~\cite[IV.2.$\beta, \g$]{book})
on the ideal $ \cL^{(1, \infty)} (\fH)$.
The local nature of the above formula stems from the
fact that the Dixmier trace vanishes on the subideal
\begin{equation*}
 \cL_0^{(1, \infty)} (\fH) \, =\, \{ T \in \cK (\fH) \, ;  \quad \sum_{i=0}^N \mu_i (T) = o (\log N) \} ,
\end{equation*}
which contains the trace class operators, and in particular
all smoothing operators on any closed manifold.
Thus,  $ \varkappa_D (a_0 , a_1, \ldots , a_p) $ depends only on the class
of the operator $a_0  [D, a_1] \cdots  [D, a_p] \, D^{-p} \in  \cL^{(1, \infty)} (\fH)$
modulo $ \cL_0^{(1, \infty)} (\fH) $, which plays the role of its symbol.

\medskip

Using the identity
\begin{equation}\label{switch1}
[D, a] \,D^{-k} \, = \, D^{-k+1} \, \big(D^k \, a \, D^{-k} \, - \, D^{k-1} \, a \, D^{-k +1}\big) \, ,
\qquad  \forall \, a \in \cA \, ,
\end{equation}
one can successively move $\,D^{-p} \,$ to the left and rewrite the 
cocycle $\varkappa_D$ in the form
\begin{equation*}
\varkappa_D (a_0 , a_1, \ldots , a_p) =
 \Tr_\om \big(\gamma a_0  (D  a_1 D^{-1}  -  a_1)  \cdots 
(D^p  a_p  D^{-p}  - D^{p-1} a_p D^{-p +1})\big) .
\end{equation*}
\medskip

In the twisted case, taking a clue from \eqref{clue},  one is led to make
the formal substitution
\begin{equation}\label{subst}
D^k \, a \, D^{-k}  \quad \longmapsto \quad
D^k \, \sigma^{-k} (a) \, D^{-k}  \, , \qquad  \forall \, a \in \cA \, ,
\end{equation}
and use the twisted version of \eqref{switch1}, namely
\begin{equation}\label{switch2}
[D, \sigma^{-k} (a)]_\s \,D^{-k} \, =
\, D^{-k+1} \, \big(D^k \, \sigma^{-k} (a) \, D^{-k} \, - \, D^{k-1} \, \sigma^{-k+1} (a) \, D^{-k +1}\big) \, ,
\end{equation}
to reverse the process of distributing $\,D^{-p} \,$ among the factors. Assuming
that the domain condition which permits the above operation is fulfilled, one
 thus arrives at the expression
  \begin{equation}\label{tHcocycle}
\varkappa_{D, \sigma} (a_0 , a_1, \ldots , a_p)\, :=
\, \Tr_\om \big(\gamma a_0  [D, \s^{-1}(a_1)]_\s \cdots  [D, \s^{-p}(a_p)]_\s D^{-p}\big) \, .
\end{equation}

This was shown in~\cite{C-M4} to be indeed 
a Hochschild $p$-cocycle, and it will be useful to
reproduce the elementary calculation that validates
this statement.
It relies on two basic properties of twisted spectral triples. The first is the obvious fact
that the $\s$-bracket with $D$ satisfies the twisted derivation rule
 \begin{equation}\label{bim2}
[D, ab]_\s \, =\,  \s(a) \, [D, b]_\s \, +\, [D, a]_\s\, b, \qquad \ a, b \in \cA .
\end{equation}
The second is the observation that
the positive linear functional on $\cA$,
\begin{equation*}
 a \, \mapsto \, \Tr_\omega (a\,|D|^{-p}) ,
\end{equation*}
is a $\sigma^{-p}$-trace on $\cA$, \ie satisfies 
\begin{equation*}\label{dixsig1}
\Tr_\omega (a\,b \, |D|^{-p}) \, =\, \Tr_\om (b \, \sigma^{-p} (a)\, |D|^{-p}), \qquad \forall \, a, b \in \cA \, .
\end{equation*} 
It is in fact a $\sigma^{-p}$-hypertrace, since  
for any $T \in \cL (\fH)$ one still has
\begin{equation}\label{dixsig2}
Tr_\omega (T \, \sigma^{-p} (a) \,|D|^{-p}) \, = \, Tr_\omega (a\, T \, |D|^{-p} ) \, .
\end{equation} 

Making use of the Leibniz rule \eqref{bim2}, one computes the Hochschild coboundary of
 $\varkappa_{D, \sigma}  \in Z^p (\cA, \cA^\ast) $ as follows:
\begin{eqnarray*}
&&b \varkappa_{D, \sigma} (a_0 , a_1, ..., a_{p+1}) \, =\, \\
&&\quad \quad = \sum_{i=0}^p (-1)^i\, \varkappa_{D, \sigma}  (a_0,..., a_i a_{i+1},..., a_{p+1})
\, + \, (-1)^{p+1} \varkappa_{D, \sigma}  (a_{p+1} a_0 , a_1, ..., a_p) \\
&&\quad \quad = \Tr_\om \big( \gamma\, a_0 \, a_1\, [D, \sigma^{-1} (a_2)]_\s \cdots \,
 [D, \sigma^{-p}  (a_{p+1})]_\s\,D^{-p}\big) \\
 &&\qquad \qquad \qquad -  \Tr_\om \big( \gamma \, a_0 a_1 \,
 [D, \sigma^{-1} (a_2)]_\s \cdots \,    [D, \sigma^{-p}  (a_{p+1})]_\s\,D^{-p} \big)\\
&&\qquad  \qquad \qquad - \Tr_\om \big( \gamma \, a_0\, [D, \sigma^{-1} (a_1)]_\s \, 
\sigma^{-1} (a_2) \cdots \,
[D, \sigma^{-p}  (a_{p+1})]_\s\,D^{-p} \big) + \ldots \\
&&\ldots +(-1)^p\, \Tr_\om \big( \gamma a_0 [D, \sigma^{-1} (a_1)]_\s \cdots \\
&& \qquad \qquad \cdots [D,  \sigma^{-p+1} (a_{p-1})]_\s \sigma^{-p+1}  (a_p) 
 [D, \sigma^{-p} (a_{p+1})]_\s D^{-p}\big) + \\
&&+ (-1)^p \, \Tr_\om \big( \gamma a_0 [D, \sigma^{-1} (a_1)]_\s \cdots \\
 && \qquad \qquad \cdots
  [D, \sigma^{-p+1} (a_{p-1})]_\s [D,  \sigma^{-p}(a_p)]_\s  \sigma^{-p}(a_{p+1}) D^{-p}\big) +\\
&& + (-1)^{p+1}\, \Tr_\om \big( \gamma  a_{p+1}  a_0\,
[D, \sigma^{-1} (a_1)]_\s  \cdots [D, \sigma^{-p}(a_p)]_\s\,D^{-p} \big) \qquad= \qquad 0 \, .
\end{eqnarray*}
The end result is $0$ because the successive terms cancel in pairs, with
the last two terms canceling each other thanks to the enhanced $\sigma^{-p}$-trace
property Eq. \eqref{dixsig2}.

\section{Ansatz for a twisted local index formula} \label{locans}

\subsection{Local index formula for spectral triples}

The local index formula that delivers in full the  Connes-Chern character 
in cyclic cohomology was developed in \cite{C-M2}, in terms of 
residue functionals that generalize Wodzicki's noncommutative
residue, and
in the framework of an abstract pseudodifferential
calculus  which gives a precise meaning
to the notion of symbol. We briefly recall the salient notions.
\medskip

Let  $(\cA , \fH , D)$ be a $(p, \infty)$-summable spectral triple
(with $D$ invertible), which satisfies the smoothness condition \eqref{R}.
We denote by $\cB$ the algebra generated by 
$\displaystyle\, \bigcup_{k \geq 0} \d^k (\cA + [D, \cA])$, and
also set
$\displaystyle\, \fH^\infty := \, \bigcap_{k \geq 0} \Dom (|D|^k)$.
We now consider linear operators
$P : \fH^\infty \ra \fH^\infty$ that admit an expansion of the form 
\begin{equation} \label{pseudodiff}
 P \, \sim \, \sum_{k \geq 0} b_{k} |D|^{s-k} \, ,  \qquad \text{with} \qquad
 b_{k} \in \cB  \, , \quad s \in \C \, ,
 \end{equation}
in the sense that 
\begin{eqnarray*} 
&& P \, - \sum_{0 \leq k <N} b_k |D|^{s-k} \,  \in \,\cB \cdot {\rm OP}^{\Re s -N}  \, ,\qquad \fl \, N > 0 \, ,\\
 \text{where} \qquad \qquad && R \in  {\rm OP}^r \quad \Longleftrightarrow \quad
  |D|^{-r} R \, \in \,  \bigcap_{k > 0} \Dom (\d^k) , \\
  \text{and} \qquad \qquad && \cB \cdot {\rm OP}^r \, := 
 \left \{ \sum_j b_j R_j \, ; \quad b_j \in \cB , \quad R_j \in  {\rm OP}^r \right \} .
  \end{eqnarray*}
Thanks to the key commutation relation (see \cite[Appendix B, Thm. B.1]{C-M2})
\begin{equation} \label{FR1}
 |D|^{s} \, b  \, \sim \,   \sum_{k \geq 0} \frac{s (s-1) \cdots (s - k +1)}{k!} \,
 \, \d^k(b) \, |D|^{s-k} \, , \qquad \fl \, b \in \cB \, ,  \quad s \in \C \, ,
 \end{equation}
these operators form a filtered algebra. 

Let now $\cD(\cA , \fH , D)$ be the algebra generated by 
$\displaystyle\, \bigcup_{k \geq 0} \nb^k (\cA + [D, \cA])$, where
$\nb$ denotes the derivation $\nb (T) = [D^2 , T]$. Its elements
play the role of {\it differential operators}. They too have a natural
order, determined by total power of $\nb$ involved in each monomial. 
Furthermore, the analogue of \eqref{FR1} holds:  \quad for any $q$-th
operator $T \in\cD^q(\cA , \fH , D)$
and $N >0$,
\begin{equation} \label{FR2}
 D^{2s} \, T  \, - \,   \sum_{0 \leq k <N} \frac{s (s-1) \cdots (s - k +1)}{k!} \,
 \, \nb^k(T) \, D^{2(s-k)} \,  \in \, {\rm OP}^{2\Re s + q-N} \, .
 \end{equation}
The intrinsic pseudodifferential calculus for the spectral triple is based on 
the algebra $\Psi^{\bullet}(\cA , \fH , D)$, generated by the
operators defined by \eqref{pseudodiff} together
with the differential operators. It is a filtered
algebra, and its quotient modulo the ideal of smoothing operators 
$\quad \Psi^{-\infty}(\cA , \fH , D)\,  = \, \bigcap_{N \geq 0} \cB \cdot {\rm OP}^{-N} $
gives the corresponding algebra of {\it complete symbols} 
$\quad \sC\sS^{\bullet}(\cA , \fH , D) \, := \, \Psi^{\bullet}(\cA , \fH , D) / \Psi^{-\infty}(\cA , \fH , D)$.
 
\medskip

Underlying the setup for  the local index formula is
the essential assumption that the spectral triple admits
a {\it discrete dimension spectrum}, to which all singularities of
zeta functions associated to elements of $\cB$
are confined; the 
postulated spectrum is a {\it discrete subset}
$\Sigma  \sbs \C$, such that the holomorphic functions
\begin{equation} \label{zeta}
  \zeta_b (z) \, = \,  \Tr (b \, |D|^{-z} ) \, , \qquad  \Re z > p \, , \quad b \in \cB ,
\end{equation}
{\it admit holomorphic extensions to $\, \C \setminus \Sigma$}. This requirement
is supplemented by a technical condition stipulating that {\it the functions   
$\, \G (z) \,  \zeta_b (z)$  decay rapidly on finite vertical strips}.
\medskip

For the sake of convenience, we shall make the stronger
assumption that the dimension spectrum is {\it simple, \ie
$\Sigma $ consists of simple poles}.  Then,
 for any $P \in \Psi^N (\cA , \fH , D) $ the zeta function
\begin{equation} \label{zeta0}
  \zeta_P (z) \, = \,  \Tr (P \, |D|^{-z} ) \, , \qquad  \Re z > p  + N\, ,
\end{equation}
can be meromorphically continued to $\C$, with simple poles in $\Sigma + N$.
Furthermore, the {\it residue functional}
\begin{equation} \label{res}
\cutint P \,: = \, {\rm Res}_{z=0} \,  \zeta_P (2z)\, , \qquad P \in \Psi^{\bullet} (\cA , \fH , D)
   \end{equation}
is an (algebraic) trace. By its very construction it
vanishes on $\Psi^{\bullet} (\cA , \fH , D) \cap \cL^1 (\fH)$, and in particular it
descends to a trace on $\, \sC\sS^{\bullet}(\cA , \fH , D)$.
\medskip

The local index formula expresses the Connes-Chern character  
$Ch^* (\cA , \fH , D) \in HC^* (\cA)$ in terms of a cocycle
in the bicomplex $\{CC (\cA), b , B\}$, 
 whose components are defined by means of the symbolic trace.
In the (invertible) odd case, these components are as follows:
for $\, q = 2\ell+1$, $\quad \ell \in \Z^+ $ ,
\begin{equation} \label{ind-odd}
\tau_{\rm odd}^q (a_0 ,\ldots ,a_q) \, = \, \sqrt{2i} \, \sum_{\bf k}  c_{q, {\bf k}} \,
 \cutint  \,a_0\,  [D, a_1]^{(k_1)}  \ldots [D, a_q]^{(k_q)} \,  |D|^{-2\vert {\bf k}\vert-q} \, ,
\end{equation}
where
\begin{equation*} \label{nabla} 
P^{(k)} \, = \, \nb^k (P)  \, ,  \qquad \fl \, P \in \Psi^{\bullet} (\cA , \fH , D) \, ,
\end{equation*}
and the coefficients are given by
\begin{eqnarray} \label{coeff}
 c_{q, {\bf k} }&=&
\frac{(-1)^{\vert {\bf k}  \vert} } {k_1 ! \ldots k_q! \, 
(k_1 +1) \ldots (k_1 + \ldots + k_q
+q)} \,  \G \left( \vert {\bf k} \vert + \frac{q}{2} \right) \, ,  \\ \nonumber \\ \nonumber
\text{where} &\ & {\bf k} \, =\, (k_1 , \ldots , k_q) \, , \qquad  \vert {\bf k}  \vert = k_1 +\ldots + k_q  \, . 
\end{eqnarray}
In the even case, with $\, q = 2\ell $, $\quad \ell \in \Z^+ $ ,
 \begin{eqnarray} \label{ind-even}
 \tau_{\rm ev}^0 (a) &=&\,  {\rm Res}_{z=0} \big(\G (z) \, \Tr (a \, |D|^{-2z})\big)  \, ,\\ \nonumber \\ \nonumber
\tau_{\rm ev}^q (a_0 ,\ldots ,a_q) &=&\,  \sum_{\bf k}  c_{q, {\bf k}} \,
 \cutint \, \g \,a_0\,  [D, a_1]^{(k_1)}  \ldots [D, a_q]^{(k_q)} \,  |D|^{-2\vert {\bf k}\vert-q} \, .
\end{eqnarray} 
 Each component  $\,\tau^q = \tau_{\rm ev / odd}^q $
 has finitely many nonzero summands, and
$\, \tau^q   \, \equiv \, 0$ for any  $\, q > p$.

\medskip

 Since the expressions  $\, \tau^q (a_0 ,\ldots ,a_q) $ are unaffected
 by the scaling $D \mapsto tD$, $ t \in \R$, we can write them in terms
 of the scale-invariant operators  
 $$  \a^k (a) \, := \, D^k a D^{-k} \, , \qquad  a \in \cA .
 $$
 Indeed, using the obvious identity
\begin{equation*} \label{recurs}
 [D, a]^{(k+1)} \, D^{-2k-3}  \, = \, D^2\left( [D, a]^{(k)} \, D^{-2k-1} \right) D^{-2}  \, - \, 
  [D, a]^{(k)} \, D^{-2k-1} \, ,
 \end{equation*}
 one verifies by induction that for any $k \geq 0$, 
  \begin{equation*} 
 [D, a]^{(k)} \, D^{-2k-1}  \, = \, \sum_{j=0}^k (-1)^j \, \binom{k}{j} \, 
 \left(\a^{2(k-j) + 1} (a)\, - \, \a^{2(k-j)} (a) \right) .
\end{equation*}
Therefore, for any $\ell \in \Z$,
 \begin{eqnarray*}
 [D, a]^{(k)} \, D^{-\ell} &=& \sum_{j=0}^k (-1)^j \, \binom{k}{j} \, 
 \left(\a^{2(k-j) + 1} (a)\, - \, \a^{2(k-j)} (a) \right) \,
 D^{2k +1 - \ell} \\ \\
 &=&  D^{2k +1 - \ell}  \, \sum_{j=0}^k (-1)^j \, \binom{k}{j} \, 
\left(\a^{\ell - 2j} (a)\, - \, \a^{\ell - 2j -1} (a) \right) \, .
\end{eqnarray*}
With the abbreviated notation
\begin{equation} \label{Sigma}
\Sigma^{(k, \ell)} (a) \, := \,\sum_{j=0}^k (-1)^j \, \binom{k}{j} \, 
 \left(\a^{\ell - 2j} (a)\, - \, \a^{\ell - 2j -1} (a) \right) \, ,
\end{equation}
 the above equality takes the form
 \begin{equation}  \label{deriv2}
 [D, a]^{(k)} \, D^{-\ell} \, = \,  D^{2k +1 - \ell}  \, \Sigma^{(k, \ell)} (a) .
\end{equation}
Successive application of the identity \eqref{deriv2} brings the
components $ \tau^q$ with $q > 0$ to the form
\begin{equation} \label{index2}
\tau^q (a_0 ,\ldots ,a_q) = \sum_{\bf k} c_{q, {\bf k} } \,
 \cutint  \, \g \, a_0\,  \Sigma^{(k_1, \, 2k_1 +1)} (a_1)   \ldots 
 \Sigma^{(k_q , \, 2(k_1 + \ldots + k_q) +q)} (a_q) \, .
\end{equation}
This formula covers the case of either parity, provided that in the odd case
we define $\g := F = D |D|^{-1}$, and incorporate the factor $\sqrt{2i}$ in the
expression \eqref{coeff} of the coefficients $ c_{q, {\bf k} }$ for $q$ odd.
\medskip

\subsection{Ansatz for the twisted case} \label{ansatz}

 Assume now that $(\cA , \fH, D, \s )$ is a $(p, \infty)$-summable
(invertible) twisted $\s$-spectral triple. 
The twisted analogue of the usual bimodule  of {\it gauge potentials} 
or {\it noncommutative differential forms} is
obviously the linear subspace $\Omega_{D , \s}^1 (\cA) \subset \cL (\fH)$ 
 consisting of operators of the form
\begin{equation*}
A = \sum_i  a_i \big(D\,b_i-\sigma(b_i)\,D\big)\, , \qquad  a_i , b_i \in \cA \, ,
\end{equation*}
which is a bimodule for the action
\begin{equation*}
a \cdot\omega\cdot b=\,\sigma(a)\,\omega\,b  \qqq \, a,\,b \in \cA \qqq \,
\omega \in \Omega_{D , \s}^1 (\cA) \, .
\end{equation*}
In the presence of the Lipschitz regularity axiom \eqref{Lreg}, one can similarly
define a bimodule $|\Om|_{D , \s}^1 (\cA) \subset \cL (\fH)$, by simply replacing $D$
with $|D|$.
Furthermore, as noted before \cf \eqref{bim2}, the map
\begin{equation*}
a \mapsto d_\sigma (a)=\,D\,a-\sigma(a)\,D
\end{equation*}
is a $\s$-derivation of $\cA$ with values in $\Omega_{D , \s}^1 (\cA) $,
and clearly, so is the map
\begin{equation*}
a \mapsto \d_\sigma (a)=\,|D|\,a-\sigma(a)\,|D| \, .
\end{equation*}
However, in order for the analogue of
the smoothness condition \eqref{R} to make sense, one needs to postulate
the existence of an extension of the automorphism $\s \in \Aut (\cA)$, and
consequently of the  $\s$-derivation $\d_\s$, to a
larger subalgebra of $\cL (\fH)$, which should contain $\Omega_{D , \s}^1 (\cA)$
as well as its higher  $\d_\s$-iterations.
\medskip

Thus, the formulation of a twisted version for the
pseudodifferential calculus is {\it not} canonical. 
Ignoring this aspect for now, let us pretend that an adequate analogue
$\Psi^{\bullet}(\cA , \fH , D , \s)$
of the algebra of pseudodifferential operators has already been constructed,
and assume that the twisted $\s$-spectral triple  $(\cA , \fH, D, \s )$ admits
a simple  discrete dimension spectrum $\Sigma \sbs \C$. 
We can then
focus on finding an  appropriate candidate for the local character cocycle.
Denoting
\begin{equation} \label{tconj}
\a_\s^k (a) \, := \, D^k \, \s^{-k} (a) \, D^{-k} \, , \qquad  a \in \cA , 
\end{equation}
and 
\begin{equation} \label{tSigma}
\Sigma_\s^{(k, \ell)} (a) \, := \,\sum_{j=0}^k (-1)^j \, \binom{k}{j} \, 
 \left(\a_\s^{\ell - 2j} (a)\, - \, \a_\s^{\ell - 2j -1} (a) \right) \, ,
\end{equation}
the analogues of the summands in Eq. \eqref{index2} are
the `residue integrals'
\begin{equation} \label{tindex2}
 \cutint  \g \, a_0\,  \Sigma_\s^{(k_1, \, 2k_1 +1)} (a_1)   \ldots 
 \Sigma_\s^{(k_q , \, 2(k_1 + \ldots + k_q) +q)} (a_q) \, .
\end{equation}
\medskip

We next define the  twisted version of the higher commutators as follows:
  \begin{equation} \label{thighcomm2}
 (a)_\s^{(k)} := \sum_{j=0}^k (-1)^j \binom{k}{j} 
 D^{2(k-j)} \s^{2j} (a)D^{2j} \, ,
 \end{equation}
respectively
 \begin{equation} \label{thighcomm}
 [D, a]_\s^{(k)} := \sum_{j=0}^k (-1)^j \binom{k}{j} 
 \left(D^{2(k-j)+1} \s^{2j} (a)D^{2j}  - 
 D^{2(k-j)} \s^{2j+1} (a)D^{2j +1}\right) \, .
 \end{equation}
 
\begin{rem} {\em Noting that}
 \begin{eqnarray*}
&& (a)_\s^{(k +1)} \, =\,  \, D^2 \, (a)_\s^{(k)} \, - \, (\s^2 (a))_\s^{(k)}\, D^2 \,  ,  \\
&& [D, a]_\s^{(k +1)}  \, =\,    \, D^2 \, [D, a]_\s^{(k)} \, - \, [D, \s^2 (a)]_\s^{(k)}\, D^2 \, ,
 \end{eqnarray*}
 {\em one could be tempted to regard the expressions \eqref{thighcomm2}, 
 \eqref{thighcomm}
 as genuine iterated twisted commutators with $D^2$. However, that would not be
 correct, because there is no guarantee that if, for instance,  $[D, a]_\s = 0$ then
 $[D, a]_\s^{(k)} = 0$ for all $k \geq 1$. A counterexample can be easily obtained
in the setting of \S 3.1, using Eq. \eqref{zcomm}. } 
 \end{rem}

 At any rate, with the above notation,
 we can now put \eqref{tindex2} in a form similar to Eq. \eqref{ind-even}.
 Indeed, the counterpart of the identity \eqref{deriv2} is
  \begin{equation*}  
 D^{2k +1 - \ell}  \, \Sigma_\s^{(k, \ell)} (a)  
 \, = \,  [D, \s^{-\ell} (a)]_\s^{(k)} \, D^{-\ell}  ,
\end{equation*}
which we then  employ to reverse the process by which the expression 
 \eqref{index2} was obtained from \eqref{ind-odd}. Keeping the same 
 notational conventions used in  \eqref{index2}, one
 thus arrives at the following
{\it Ansatz} for the twisted version of the local character cocycle:
  \begin{eqnarray} \label{tindex} 
&&\quad \tau_\s^q (a_0 ,\ldots ,a_q) \, =   \\ \nonumber
&& \sum_{\bf k}  c_{q, {\bf k} } \, 
 \cutint \, \g \, a_0\,  [D , \s^{-2k_1 - 1} (a_1)]_\s^{(k_1)}  \ldots 
 [D ,  \s^{-2(k_1 + \ldots + k_{q}) -q} (a_q)]_\s^{(k_q)} \, \vert D \vert^{-2\vert {\bf k} \vert-q} .
\end{eqnarray}  
There is an immediate obstruction for this formula to define a $(b, B)$-cocycle,
which arises from the   
$B$-coboundary of $\tau_{\rm odd}^1$. Indeed,
  \begin{eqnarray*} 
   B\tau_\s^1 (a) &=&  \,\sum_{k \geq 0}  c_{1, k} \, 
 \cutint \, F \, [D , \s^{-2k - 1} (a)]_\s^{(k)}  \, D^{-2k-1} \, = \, \sum_{k \geq 0}  c_{1, k} \, 
  \cutint \, F \, \Sigma_\s^{(k, \, 2k +1)} (a)  \\
  &=&  \sum_{k \geq 0}  c_{1, k} \, \sum_{j=0}^k (-1)^j \,   \binom{k}{j} \,  \left(\cutint \, F \,
 a_\s^{ 2(k - j) +1} (a)\, - \, \cutint \, F \,\a_\s^{2(k- j)} (a) \right) \\
  &=&  \sum_{k \geq 0}  c_{1, k} \, \sum_{j=0}^k (-1)^j \,   \binom{k}{j} \,  \left(\cutint \, F \,
 \s^{ -2(k - j) -1} (a)\, - \, \cutint \, F \,\s^{-2(k- j)} (a) \right) .
\end{eqnarray*}
This expression vanishes if 
  \begin{equation} \label{siginv}
 \cutint \, F \,\s (a) \, = \, \cutint \, F \, a   \, , \qquad \fl \, a \in \cA \, ,
\end{equation}
or equivalently
 \begin{equation*} 
  \cutint \, [D, a]_\s \, |D|^{-1} \, = \, 0 \, , \qquad \fl \, a \in \cA \, ,
\end{equation*}

\medskip

 As it will become apparent in the next section, Eq. \eqref{siginv}
 has something in common with the Selberg principle
 for orbital integrals of reductive Lie groups. 
 We shall show later that for a special class of conformally twisted
spectral triples there are no higher obstructions
to the validity of the Ansatz.
 \bigskip

\section{Conformal geometry and twisted spectral triples} 
\label{conforst}

 In order to shed some light on the nature and plausibility of the
 above setup for the Ansatz,
 we examine in this section some authentic examples of twisted
 spectral triples arising in conformal geometry.

\subsection{Transversely conformal spectral triple} \label{tcst}
Let  $M$ be a smooth connected closed spin manifold of
dimension $n$. To each riemannian metric
 $g$ on $M$ one can canonically associate a Dirac operator $D = \Dirac_g$
 acting on the Hilbert space $\fH = \fH_g : = L^2(M,\cS_g)$ of $L^2$-sections of
 the spin bundle $\cS= \cS_g$,
and thus a corresponding spectral triple 
$(C^\infty(M), \fH,\,  D)$ over the algebra $C^\infty(M)$.
Assume now that $M$ is endowed with a
conformal structure $[g]$, consisting of all riemannian metrics  conformally equivalent 
to a given riemannian metric $g$. Let $SCO (M, [g])$ be the 
group of diffeomorphisms
of $M$ that preserve the conformal structure, the orientation and the spin structure.
It is a Lie group, and we denote by
$G = SCO (M, [g])_0$ its connected component of the identity. We then form the
discrete crossed product algebra $\cA_G =\,C^\infty(M)\rtimes G$.
This algebra consists of
 finite sums of the form
$$ a=\,\sum_\Gamma
f_\phi\, v_\phi \,, \qquad f_\phi \in \,C^\infty(M) , \quad \phi \in G ,
$$
with the product rule determined by
\begin{equation*}
v_\phi\,f=\,(f \circ\phi^{-1})\, v_\phi \,,\quad
v_\phi \, v_\psi= v_{\phi\,\psi}\, .
\end{equation*}
 It can be represented by bounded linear operators on the 
Hilbert
space $\fH=\,L^2(M,\cS)$
of $L^2$-sections of the spin bundle $\cS$,  by letting
a function $f \in C^\infty(M)$ 
act as the multiplication operator
\begin{equation} \label{vrep1}
\pi(f)\,(u) \, =\,f\, u  \, ,\qquad 
u \in L^2(M, \cS) \, , 
\end{equation}
and the diffeomorphisms  $\phi \in G$ act
as translation operators
\begin{equation} \label{vrep2}
\pi(v_\phi)\,(u) \equiv V_\phi\,(u) \,: =\,\tilde{\phi} \circ u \circ \phi^{-1} \, ,\qquad 
u \in L^2(M, \cS) \, , 
\end{equation}
where $\tilde{\phi}$ is the canonical lift of $\phi$ to an automorphism of $\cS$; such 
a lift is well-defined, not just modulo $\Z/2\Z$, for any $\phi \in  SCO (M, [g])_0$. 
 To make $G$ act by unitary operators, 
one needs to replace each operator  $V^{-1}_\phi$, $\phi \in G$,
 by the operator
\begin{equation} \label{urep}
U^{-1}_\phi\,(u) \, = \,  e^{-n h_\phi} \, V^{-1}_\phi \,(u) \,
=\, e^{-n h_\phi} \, \tilde{\phi}^{-1} \circ u \circ \phi \, ,\qquad 
u \in L^2 (M, \cS) \, ,  
\end{equation}
where $\, h_\phi \in C^\infty(M)$ is determined by the conformal factor via the equation
\begin{equation} \label{hphi}
\phi^* (g) = e^{-4 h_\phi}\, g\, .
\end{equation}
Indeed, using the fact that the riemannian volume forms are related by the equality 
\begin{equation} \label{volchange}
 \vol_{\phi^\ast(g)} = e^{-2n h_\phi}\,\vol_g \, ,
\end{equation}
and denoting the fiberwise norm by
$| \cdot |$, one easily checks that $U^{-1}_\phi$ is unitary:
 \begin{eqnarray*}  
 ||U^{-1}_\phi\,(u)||^2 &=& \int_M e^{-2n h_\phi} \, 
 |\tilde{\phi}^{-1} (u \circ \phi) |^2 \, \vol_g \, = \, \int_M |\tilde{\phi}^{-1} (u \circ \phi)|^2 \,
 \phi^*(\vol_g ) \\ 
 &=&\,  \int_M | u |^2 \,\vol_g \, = \, || u ||^2  , \, \qquad \fl \, u \in L^2 (M, \cS) \, .
\end{eqnarray*}

 \begin{lem}  For any $\, \phi \in  G = SCO (M, [g])_0$ and with $\, D =  \Dirac_g$,
 one has
 \begin{equation} \label{Dircconfdiff3}
U^*_\phi \circ D \circ U_\phi  \,
=\, e^{h_\phi}\circ D  \circ e^{h_\phi}   \, .
  \end{equation}
\end{lem}

\proof   Via the natural identification $\, \b^{\phi^\ast(g)}_g$
corresponding to the change of metric,  defined as in \cite{bourg},
the Dirac operator $\Dirac_{\,\phi^\ast(g)} $ can be implemented as an operator
\begin{equation*} \label{Dirchange}
D_{\phi^\ast(g), \, g}\, = \, 
\left(\b^{\phi^\ast(g)}_g \right)^{-1} \circ \Dirac_{\,\phi^\ast(g)} \circ \b^{\phi^\ast(g)}_g
\end{equation*}
acting on the sections of the bundle $\, \cS_g$. It is explicitly given by the formula 
\begin{equation} \label{Dirconf}
D_{\phi^\ast(g), \, g} \, =\, e^{(n+1)h}\,\Dirac_g  \, e^{(-n+1)h} \, .
\end{equation}
On the other hand, as differential operator,
 \begin{equation} \label{Dirchangediff}
 D_{\phi^\ast(g), \, g}\, = \,  V^{-1}_\phi \circ \Dirac_g \circ V_\phi   \, .
  \end{equation}   
Combining \eqref{Dirconf} and \eqref{Dirchangediff} one obtains
 \begin{equation*} \label{Dircconfdiff1}
  V^{-1}_\phi \circ \Dirac_g \circ V_\phi  \, =\, e^{(n+1)h}\,\Dirac_g  \, e^{(-n+1)h}   \, ,
  \end{equation*}   
or equivalently
 \begin{equation*} \label{Dircconfdiff2}
e^{-n h}\circ  V^{-1}_\phi \circ \Dirac_g \circ V_\phi  \circ e^{nh}\,
=\, e^{h}\circ \Dirac_g  \circ e^{h}   \, .
  \end{equation*}
\endproof

Let $\sigma$ be the algebra automorphism of
$ \cA_G$ defined on generators by
\begin{equation} \label{defsigma}
 \s(f\, v^{-1}_\phi)=\,e^{-2 h_\phi}\,f\,v^{-1}_\phi
  \,, \qquad f \in \,C^\infty(M) , \quad \phi \in G .
\end{equation}
\medskip
 
 \begin{lem} \label{bdcom1} 
The twisted commutators
\begin{equation*} \label{ucomm}
[D, \pi(a)]_\s := D \circ \pi(a) \,-\, \pi(\s(a))\circ D ,
\qquad a \in \cA_G \, ,
  \end{equation*} 
  are bounded.
\end{lem}
\smallskip

\proof It suffices to check the claimed property for  
 $a\, =\, e^{-nh_\phi }\, v^{-1}_\phi$. In that case one has  
\begin{equation*} 
[D, \pi(a)]_\s \, = \, D \circ U_\phi^\ast \,  - \,e^{-2 h_\phi} \circ U_\phi^\ast \circ D
\, = \, \left(D\,  - \,e^{-2 h_\phi} \circ U_\phi^\ast
  \circ D \circ U_\phi\right) \circ U_\phi^\ast  \, .
 \end{equation*}
 In view of Eq. \eqref{Dircconfdiff3}, it follows that
\begin{eqnarray}  \label{DUphis}
\qquad [D, U_\phi^\ast ]_\s &=&  
\left(D\,  - \, e^{- h_\phi  } \circ
 D \circ e^{h_\phi }\right) \, U_\phi^\ast \, = \,
 = \, - \, e^{- h_\phi} \, [D,  e^{h_\phi}] \circ U_\phi^\ast \\ \notag
 &=&
 - c(dh_\phi) \, U_\phi^\ast .
 \end{eqnarray}
\endproof

 For further reference, we note that as a consequence of \eqref{DUphis} one has
\begin{equation}  \label{expcomm} 
[D, f\, U_\phi^\ast ]_\s \, = \, c(df \, - \, f \,dh_\phi) \, U_\phi^\ast   \, ,
 \end{equation} 
and in particular
 \begin{equation}  \label{zcomm} 
[D, e^{h_\phi }\, U_\phi^\ast ]_\s \, = \, 0 \, .
 \end{equation}

\begin{prop}\label{tst}
The algebra $\cA_G =\,C^\infty(M)\rtimes G$, endowed with the 
automorphism $\s \in \Aut \cA_G$ and the representation $\pi$ on the Hilbert
space $\fH=\,L^2(M,\cS)$,  
together with the Dirac operator $D = \Dirac_g$, define an
$(n, \infty)$-summable $\s$-spectral
triple $(\cA_G , \fH, D, \s)$,
which moreover satisfies the strong
 Lipschitz-regularity property
 \begin{equation}  \label{SLR}
|D|^{-t} \left(|D|^t \, a\, - \,\s^t (a)  \,|D|^t \right)\, \in \, \cL^{(n,\infty)}(\fH)\, ,
\qquad \fl \, t \in \R .
\end{equation}
\end{prop}
\smallskip
 
 \proof  The boundedness property \eqref{boundedcom} was verified in Lemma \ref{bdcom1}.
To prove Lipschitz regularity, one notes that, when viewed as
pseudodifferential operator, $U^*_\phi \circ D \circ U_\phi$ has
principal symbol
\begin{equation} \label{prs1}
\prs (U^*_\phi \circ D \circ U_\phi)(x, \xi) \, = \, e^{2 h_\phi} c (\xi) \, , \quad \xi \in T^*_x M
\end{equation}
where $c (\xi) \in \End (\cS_x)$ stands for the Clifford multiplication by $\xi$. 
Furthermore, for any $t \in \R$, the principal symbol of
$U^*_\phi \circ |D|^t \circ U_\phi$ is  
\begin{equation} \label{prs2}
\prs (U^*_\phi \circ|D|^t \circ U_\phi)(x, \xi) \, = \, e^{2 t h_\phi} ||\xi||^t \, ,
\end{equation}
since 
  \begin{equation*}  
U^*_\phi \circ|D|^t \circ U_\phi   \, = \, |U^*_\phi \circ D \circ U_\phi|^t  \, .
 \end{equation*}
Now 
\begin{equation*}  
|D|^t \circ  U^*_\phi \, - \,\s^t (U^*_\phi)  \circ |D|^t \, = \,\left( |D|^t  \, - \,
e^{-2t h_\phi}\,U^*_\phi \circ |D|^t  \circ U_\phi \right) \circ  U^*_\phi \, ,
\end{equation*}
and by Eq. \eqref{prs2},
\begin{equation*}  
 \prs \left( |D|^t  \, - \,
e^{-2 t h_\phi}\,U^*_\phi \circ |D|^t  \circ U_\phi \right)  \, = \, 
||\xi||^t  \, - \,
e^{-2 t h_\phi}\, e^{2 t h_\phi} ||\xi||^t  \, = \, 0 .
\end{equation*}
Thus,  the operator $\, |D|^t  \, - \,
e^{-2 t h_\phi}\,U^*_\phi \circ |D|^t  \circ U_\phi \,$ is
pseudodifferential of order $t - 1$, hence 
its product by $|D|^{-t}$ is of order $-1$ and therefore in $\cL^{(n ,\infty)}$.
 \endproof
\medskip

\begin{rem} \label{flatst}
In the same fashion, 
$\R^n$ with its standard metric $g_0$,
together with the flat Dirac operator $D_0 = \Dirac_{g_0}$, give rise to
 the $(n, \infty)$-summable non-unital $\s$-spectral
triple $(\cA_{G_0} , \fH_0, D_0, \s)$, where $G_0 = CO (\R^n, g_0)$, and
$\cA_{G_0} = C_c^\infty (\R^n) \rtimes G_0$. 
\end{rem}  
\medskip

According to the Ferrand-Obata theorem (\cf ~\cite{fer} for a complete proof),
the conformal group $CO (M, [g])$ of a (not necessarily closed) manifold $M$ of
dimension $n \geq 2$ is {\em inessential}, \ie reduces to
the group of isometries for a metric in the conformal class $[g]$,
except when $M^n$ is conformally equivalent
 to the standard sphere $S^n$ or to the standard Euclidean space $\R^n$.
Correspondingly, the only twisted spectral triples arising from the
above construction which are not isomorphic to ordinary spectral triples
are those associated to the $n$-sphere and to the flat $n$-space.  

\bigskip

\subsection{Transverse noncommutative residue} \label{confind}
 With the same assumptions as in the preceding subsection, let 
 $\Psi^{\bullet} (M ; \cS)$ denote 
 the algebra of classical pseuododifferential operators
acting on the sections of the spin bundle. 
The group $G$ acts on $\Psi^{\bullet} (M ; \cS)$ in the natural fashion:
\begin{equation} \label{Gpsi} 
 \phi \cdot P \, := \,  V_\phi \, P \, V^{-1}_\phi , \qquad    \phi \in G , \quad 
 P \in \Psi^{\bullet} (M ; \cS).
  \end{equation}
One can thus form the crossed product algebra $\Psi^{\bullet} (M ; \cS) \rtimes G$.
The representation $\pi : \cA_G \ra \cL (\fH)$ extends in a tautological manner
 to a representation of the enlarged algebra
$\Psi^{\bullet} (M ; \cS) \rtimes G$ by densely defined linear operators on $\fH = L^2 (M, \cS)$,
which will still be denoted by $\pi$:
\begin{equation} \label{vrep3}
\pi(P \, v_\phi)\,(u) \,: =\,P \big(V_\phi(u)\big) \, = \,P (\tilde{\phi} \circ u \circ \phi^{-1}) \, ,\qquad 
u \in C^\infty(M, \cS) \, .
\end{equation}

We set out to show that
given any $\, \cP \in \Psi^N (\cA_G)$ the zeta function
 \begin{equation} \label{zeta1} 
  \z_{\cP} (z) \, := \, \Tr (\cP \, |D|^{-z} )  \, , \qquad \Re z  > n + N
  \end{equation}
 can be meromorphically continued to the whole complex plane;
  by linearity, it suffices to take $\cP = P\, V_\phi $, with
  $P \in \Psi^N (M ; \cS)$ and $\phi \in G$.
  \medskip
  
 If  $G$ is inessential, $\phi$ is an isometry and the statement can 
 be proved via the Mellin transform and  heat kernel
 asymptotics (see \cite[\S 6.3]{B-G-V}). A more direct proof, given in~\cite{Dave},
 relies on the stationary phase method 
 (\cf \eg \cite[Thm. 7.7.5]{Horm}), applied to a phase function whose expression
 in local charts $U$ covering a tubular neighborhood of the fixed point set 
 $M_\phi$ is of the form
 \begin{equation} \label{phase1}
  f (x, \xi) \, = \, \langle x - \phi (x) , \xi \rangle  \, , \qquad x \in U , \quad \xi \in \R^n , \, 
  ||\xi|| = 1.
  \end{equation}
By restriction to the fibers of the
normal bundle to $M_\phi$, this function gives rise to a family of fiberwise
phase functions, each having a single non-degenerate stationary point. Using
the stationary phase for this family, it is shown in \cite[Prop. 2.4]{Dave}
that the zeta function $ \z_{V_\phi \, P}$ has a
meromorphic extension to $\C$ whose poles are at most 
simple and located at the points
$\, \,  z_k  = N + n_\phi - k$,  $\, \,  k   \in \Z^+$, \, where $\, \, n_\phi = \dim M_\phi$.

 \medskip
 
In the sphere case, by Liouville's theorem the group of conformal automorphisms
$CO(S^n , [g])$
coincides with the group  $M(n) \cong PO(n+1, 1)$ of M\"obius transformations 
in dimension $n$, and $\, G =  SCO(S^n , [g])$ is its connected component.
 Now if $\phi \in G$ is {\it elliptic}, \ie  conjugate
 to an element in the maximal compact subgroup $O(n+1)$,
by replacing $D$ in formula \eqref{zeta1} with a conjugate by
 a unitary operator we can reduce to the isometric case.
 
The non-elliptic diffeomorphisms  $\phi \in M(n)$ fall into
two classes: {\em hyperbolic} and {\em parabolic} (see \cite[\S 2]{kulk}).

A {\it hyperbolic} transformation $\phi \in M(n)$
has two distinct fixed points, say $ x^+$ and $x^- $, and its tangent map at each
of these points  $d \phi_{x^{\pm}} : T_{x^{\pm}} S^n \rightarrow T_{x^{\pm}} S^n$
is represented by an element of $O (n) \times \R^+$,
with multiplier $\mu^{\pm}$, with $\mu > 1$. Because of the nontrivial
multiplier, the phase function \eqref{phase1} 
has no critical points away from the zero section of the cotagent
bundle $ T^{\ast}  S^n $.  
The stationary phase principle, in its most basic form (\cf \cite[Thm. 7.7.1]{Horm}) and
utilized in the same manner 
as in~\cite{Dave}, implies then that
the zeta function $ \z_{V_\phi \, P}$ extends to an entire function.
\smallskip

A {\it parabolic} transformation $\phi \in M(n)$ has a single fixed point $x_0 \in S^n$, 
however
$\det (\id - d \phi_{x_0}) = 0$. Accordingly,
the phase function has $0$ as
its only critical value, and the corresponding critical set is
  \begin{equation} \label{crit1}
 C_\phi \, = \, \{ (x_0, \xi) \,   \vert \,  \xi \in \R^k , \quad d\phi_{x_0} (\xi)=\xi \} .
  \end{equation}
The generalized stationary phase gives
 an asymptotic expansion 
 \begin{equation} \label{asymp0}
 \int_{ ||\xi|| = 1} e^{ i r \langle x - \phi (x) , \, \xi \rangle} \, a (x, \xi) \, d^{n-1}\xi \, d^n x\, \sim \,
 \sum_{\alpha} \sum_{j=0}^\infty \sum_{k=0}^{2n - 1}\, \d_{ j, k} (a) \, 
 r^{\alpha - j} \, \log^k r \, ,
 \end{equation}where $\alpha = \frac{1}{2}\dim C_\phi -n$   and
the distributions $\d_{ j, k}$ are supported in $ C_\phi$. 
Following the same line of arguments as in \cite{Dave}, but using stereographic
coordinates instead of normal coordinates,
one obtains the desired
meromorphic continuation of the zeta function $ \z_{V_\phi \, P}$, with (at most)
simple poles at the points   
$$ z_k  = N + \frac{1}{2}\dim C_\phi - k , \qquad \text{where} \qquad k   \in \Z^+ .
$$
\medskip

We summarize the conclusion of the preceding discussion as follows.
 \medskip
 
  \begin{thm} \label{ibp0} 
  For any $P \in \Psi^N (M^n ; \cS) $ and any $\phi \in G $, 
  the associated zeta function  $ \z_{P\, V_\phi}$ has a meromorphic
 extension to $\C$. Moreover,
    \begin{itemize}
      
  \item[$1^0$]  if $\phi \in G$ is elliptic,  
  then the poles of  $ \z_{V_\phi \, P}$ are at most simple and
  are   located at the points
  $\, z_k  = N + \dim M_\phi - k$, $\quad k   \in \Z^+ $;
    
  \item[$2^0$] if $\phi \in G$ is hyperbolic, then  $ \z_{V_\phi \, P}$
  is entire;
  
  \item[$3^0$]  if $\phi \in G$ is parabolic, 
  then the poles of  $ \z_{V_\phi \, P}$  are at most simple and
  are located at the points
  $\, z_k  = N + \frac{1}{2}\dim C_\phi  - k$, \quad $ k   \in \Z^+ $.
  \end{itemize}
 \end{thm}
 \medskip
 
 This provides the {\em transverse noncommutative residue} functional
 \begin{equation*}
   \cutint  \, \cP \, = \, \Res_{z=0} \, \z_{\cP} (2z) \, , \qquad  \cP \in  \Psi^{\bullet} (M ; \cS) \rtimes G ,
 \end{equation*}
 which satisfies a property analogous to the Selberg Principle.
 
 \begin{cor} \label{selberg}
For any hyperbolic transformation $\phi \in G$ and any  $P \in \Psi^N (S^n ; \cS) $,
   \begin{equation*}  \label{hyp}
\qquad  \cutint \,  P \,V_\phi \,  = \, 0\,  .
 \qquad \qquad 
  \end{equation*} 
 \end{cor}
  \medskip

As another consequence, one can explicitly compute  
the candidate for the Hochschild character 
given by Eq. \eqref{tHcocycle}, and thus directly verify that it
gives the expected result.
 
 \begin{prop} \label{tHclass1}
   The local Hochschild cocycle of the transversely conformal 
 $\s$-spectral triple $(\cA_G , \fH_g, \Dirac_g, \s)$ associated to a closed spin manifold $M^n$
 modulo the conformal group $G = SCO(M, [g])$
 is a cyclic cocycle whose
 periodic cyclic cohomology class coincides with the
 transverse fundamental class $[M/G]$.
 \end{prop}
   \medskip
   
\proof Let $\, a^k \, = \, f_k \, U^*_{\phi_k} \in \cA_G$, $k = 0, 1, \ldots , n$. The
integrand in the formula \eqref{tHcocycle}, 
$$
a^0 \, [D , \s^{-1}(a^1)]_\s\, \cdots \,
[D , \s^{-n} (a^n)]_\s \,\vert D\vert^{-n} \, ,
$$
can be put in the form $\, P \, V_\phi $ with $\, P \in \Psi^{-n} (S^n ; \cS) $ and
$\, \phi^{-1} =   \phi_n \circ \ldots \circ \phi_0$. It follows from Prop. \ref{ibp0}, 
specialized to the case when $N = -n$, that
the zeta function $ \z_{P \, V_\phi}(z)$ has no pole at $z=0$.
Therefore  
  \begin{equation*}   
\qquad  \cutint  P \, V_\phi \,  = \, 0\,  , \qquad \text{unless} \quad \,  \phi = \id ,
 \qquad \qquad 
  \end{equation*} 
\ie the cocycle \eqref{tHcocycle} is localized at the identity.
Employing Getzler's symbol calculus for
asymptotic operators as indicated in \cite[Remark II.1]{C-M2}, and
using the expression \eqref{expcomm} 
of the twisted commutators, the Hochschild cocycle \eqref{tHcocycle} can be explicitly computed. 
The end result is a cyclic cocycle, which is easily seen to differ by a 
coboundary from the standard transverse fundamental cocycle (\comp \cite[Thm 3.11]{C-M4}) 
\begin{eqnarray*} \displaystyle
\tau_{M/G} (f_0 U^*_{\phi_0} , \ldots , f_n U^*_{\phi_n} )= \left\{ \begin{matrix}
\int_M f_0 \, d(f_1 \circ \phi_0) \wedge \ldots \wedge  d(f_n \circ \phi_{n-1} \circ \ldots \circ \phi_0) ,\\ 
\qquad  \qquad \qquad \qquad \text{if} \qquad
\qquad \phi_n \circ \ldots \circ \phi_0 = \id ;\\
\qquad \qquad \quad 0 \qquad \qquad \text{otherwise} \, .
\end{matrix} \right.
\end{eqnarray*}
\endproof
\medskip

\begin{rem} \label{inesstst}
{\em If  $M$ is closed and $G$ inessential, hence compact,
the corresponding
 $\s$-spectral triple is a conformal perturbation, \cf  \S \ref{confpert},
 of an equivariant spectral triple. Its Connes-Chern character,
 given by \eqref{chpert}, can be explicitly computed from the local
 index formula of \cite{C-M2} by employing an equivariant version of
 the Getzler symbol calculus, or the equivariant heat kernel techniques in \cite{B-G-V}.}
 \end{rem} 
\bigskip
 
\subsection{Transverse similarities} \label{confRn}
Endow $ \R^n$ with the Euclidean metric $\, g_0 = \sum_{i=1}^n d x^i \ot d x^i $.
The group $G = {\rm Sim} (n)$ of conformal (or similarity)
transformations of the Euclidean $n$-space is generated by rotations,
translations by vectors 
$\, y \in \R^n$,
$$
\tau_y (x) \, = \, x - y , \qquad \fl \, x \in \R^n ,
$$
and homotheties 
$\, \rho_\lb$, $\, \lb > 0$, 
$$
\rho_\lb (x) \, = \, \lb^{-1} x , \qquad \fl \, x \in \R^n .
$$
The only non-isometries are the homotheties,
\begin{equation*}
\rho_\lb^\ast g_0 \,=\,  \lb^{-2} g_0 \, .
\end{equation*}
\ie in the notation of \S \ref{tcst} the corresponding conformal factor is
\begin{equation*}
e^{- 4h_{\rho_\lb}}(x) \, = \,  \lb^{-2}  \, , \quad \fl \, x \in \R^n .
\end{equation*}
With $\cA_G := C_c^\infty (\R^n)  \rtimes G$, the definition \eqref{defsigma}
of the automorphism $\s \in \Aut \cA_G$ specializes to
 \begin{equation} \label{charaut}
 \s(f\, U_\phi)\, =\,\mu (\phi)\,  
 \,f\,U_\phi \, , \qquad \phi \in G \, ,
\end{equation}
where  $\, \mu : G \ra \R^+$ is the character determined by the {\it multiplier} of
the similarity transformation:
$$ \mu (\phi) =1\quad \text{if} \quad \phi \in O(n) , \quad
\mu (\tau_y) = 1, \, \fl \, y \in \R^n , \quad
\text{and} \quad \mu (\rho_\lb) \, = \, \lb , \, \fl \, \lb  > 0 .
$$
 Also, the covariace relation \eqref{Dircconfdiff3} becomes
\begin{equation} \label{Dcov}
   U^{-1}_\phi \circ D \circ U_\phi \, = \, \mu(\phi)\,D \, , \qquad \phi \in G\, .
  \end{equation}   
\medskip

Let $ \Psi_c^{\bullet} (\R^n ; \cS)$ denote the algebra of classical pseudodifferential operators
with $x$-compact support. Since the conformal factors
are constant, one can easily 
 extend $\s$ to an automorphism of  $\Psi_c^{\bullet} (\R^n ; \cS) \rtimes G$, by simply setting
\begin{equation} \label{Psiaut}
 \s(P\, U_\phi)\, =\,\mu (\phi)\,  
 \,P\,U_\phi \, , \qquad \phi \in G , \quad P \in \Psi_c^\bullet (\R^n, \cS) .
\end{equation}
Indeed,
\begin{eqnarray} \label{extends}
&&\qquad \s(P\, U_\phi) \cdot \s(Q \, U_\psi)\,=\,
\big( \mu (\phi)\,  \,P\,U_\phi \big) \cdot  \big( \mu (\psi)\,  \,Q\,U_\psi \big) \, = \\ \notag
&&\quad =\, \mu (\phi \psi) \, P \cdot  (U_\phi Q U^{-1}_\phi ) \, U_{\phi \psi}\, =\,
  \s\big(P \cdot  (U_\phi Q U^{-1}_\phi)  \, U_{\phi \psi}\big) \, = \,
 \s(P\, U_\phi \cdot  Q \, U_\psi) \,.
 \end{eqnarray} 
 \medskip

 \begin{prop} \label{trinv} The residue functional 
 $\displaystyle \, \cutint : \Psi_c^{\bullet} (\R^n ; \cS) \rtimes G \ra \C \,$ is a
 $\s$-invariant  trace.
\end{prop}

\proof  The $\s$-invariance of the residue is a consequence of 
Corollary \ref{selberg} (Selberg Principle).
Indeed, let  $ \cP= PU_\phi \in \Psi_c^{\bullet} (\R^n ; \cS) \rtimes G$.
 If $\mu(\phi) \neq 1$ 
then $\phi$ has a unique
fixed point $x_0 \in \R^n$, at which $d\phi_{x_0} = \mu(\phi) \id$. Thus, there are
no fixed points on $T^\ast \R^n$, hence the zeta function
$\z_\cP (z)$ is entire (see \S \ref{confind}).
 On the other hand,  if $\mu (\phi) = 1$ then
$\s (\cP) = \cP$.

To prove that the residue functional is a trace, let 
$ PU_\phi , \, QU_\psi \in \Psi_c^{\bullet} (\R^n ; \cS) \rtimes G$. One has
\begin{eqnarray*}
\Tr \big(PU_\phi \, Q U_\psi \, |D|^{-z}\big) &=&\Tr \big(Q U_\psi \, |D|^{-z}  P\, U_\phi \big) = \\
&=& \Tr \big(Q U_\psi \, P |D|^{-z} \, U_\phi \big) +  \Tr \big(Q U_\psi \, [|D|^{-z},  P]\, U_\phi \big)
\end{eqnarray*}
Using the
identity Eq. \eqref{FR1} to express the commutator $[|D|^{-z},  P]$, one sees that
$$
\Res_{z=0}  \Tr \big(Q U_\psi \, [|D|^{-z},  P]\, U_\phi \big) \, = 0 \, ,
$$
hence,
\begin{eqnarray*}
\cutint PU_\phi \, Q U_\psi = \Res_{z=0} \Tr \big(Q U_\psi \, P |D|^{-z} \, U_\phi \big) =
 \Res_{z=0} \Tr \big(Q U_\psi \, P U_\phi \, U_\phi^{-1} |D|^{-z} \, U_\phi \big) .
\end{eqnarray*}
By Eq. \eqref{Dcov}, $\quad  U_\phi^{-1} |D|^{-z} \, U_\phi = \mu(\phi)^{-z}  |D|^{-z}$, 
therefore
\begin{eqnarray*}
\cutint \, PU_\phi \, Q U_\psi =
 \Res_{z=0} \left(\mu(\phi)^{-z} \, \Tr \big(Q U_\psi \, P U_\phi \,  |D|^{-z} \big)\right) =
 \cutint \, Q U_\psi \,  PU_\phi .
\end{eqnarray*}
\endproof

\begin{rem} One can explicitly verify in this specific case 
that the cochain  \eqref{tindex} of the Ansatz does satisfy the cocycle identity
$$
 b\, \tau^{q-1}_\s \, + \, B\, \tau^{q+1}_\s \, = \, 0 .
 $$
  \end{rem}
Indeed the direct, albeit lengthy, computations by which the cocycle
identity is checked in the beginning
of the proof in \cite[Theorem II.1]{C-M2}  
can be reproduced
almost {\em verbatim}.  Once the commutator with $D^2$ is substituted by the
twisted commutator
\begin{equation} \label{nablas}
\nabla_\s (\cP) \,= \, D^2\, \cP \, - \, \s^2 (\cP) \, D^2  \, ,
\end{equation} 
and usual iterated Leibniz rule is replaced with
its twisted version
     \begin{eqnarray} \label{tdergen}
&&  \nabla_\s^m (\cP_1 \cP_2 \cdots \cP_q )\, =  \,
\sum_{m_1 + \ldots + m_q = m}  \frac{m!}{m_1! \cdots m_q!}  \cdot    \\ \nonumber
&&\qquad \qquad \nabla_\s^{m_1} ( \s^{2(m_2 + \ldots + m_q)}(\cP_1) ) \, 
  \nabla_\s^{m_2} (\s^{2(m_3 + \ldots + m_q)} (\cP_2) ) \cdots \,
\nabla_\s^{m_q}(\cP_q)  \, ,
 \end{eqnarray} 
the ``integration by parts''  property, which is repeatedly used in those calculations,
becomes a consequence of Proposition \ref{trinv}. As a simple illustration,
\begin{eqnarray*} 
  \cutint \, \nb_\s (\cP)  \, D^{-2} \, = \, \cutint \, D^2 \,\cP \, D^{-2} \, - \, \cutint \, \s^2 (\cP)
  \, = \, \cutint \, \cP \, - \,  \cutint \, \s^2 (\cP)  \, = \, 0 \, .
  \end{eqnarray*}
 \smallskip
 
A different approach, which provides a
complete verification of the Ansatz in greater generality,
makes the object of the section that follows.
 \bigskip
 
\section{Twisting by scaling automorphisms and the local index formula} \label{transim}

By abstracting the essential features of the preceding example,
we define in this section a general class of spectral triples
twisted by scaling automorphisms,
for which we shall prove the validity of the  Ansatz in its entirety.

\subsection{Scaling automorphisms} Motivated by the Euclidean similarities in
\S \ref{confRn}, we introduce the following abstract version of a
spectral triple twisted by similarities.

\begin{defn} \label{consim}
Let  $(\cA , \fH , D)$  be a  spectral triple over the 
(non-unital) involutive algebra $\cA$.
 A {\em scaling automorphism} of  $\, (\cA , \fH , D)$ is defined by a unitary operator
 $U \in \cU(\fH)$ such that 
  \begin{equation} \label{sim}
U\, \cA \, U^\ast \, = \,\cA \, , \quad \text{and} \quad
 U\, D \, U^\ast \,= \,\mu(U) \, D , \quad \text{with} \quad \mu(U) > 0 .
\end{equation}
\end{defn}
Scaling automorphisms form a group 
$\, \Sim(\cA , \fH , D)$, endowed by construction with a
 {\em scaling character} $\mu :  \Sim(\cA , \fH , D) \ra \R^+ $.
 Its subgroup $\Ker \mu$ consists of the 
 {\em isometries} of $\, (\cA , \fH , D)$, and we will be denoted
$\Isom(\cA , \fH , D)$.

  \medskip
 
For the clarity of the exposition it will be convenient to assume $D$ invertible.
This can always be achieved by passage to the invertible double, \cf  \S \ref{invdbl},
\begin{equation*}
\wt{D} \, = \, D \hot \id\, + \, \id \hot F_1 \, .
\end{equation*}
However, in doing so the similarity condition \eqref{sim} cannot be exactly
reproduced. Instead, it takes the modified form
\begin{equation} \label{psim}
 U\, \wt{D} \, U^\ast \, = \,\mu(U) \,  \wt{D} \hot \id\, + \,(1 - \mu(U) )\,  \id \hot F_1 \, .
\end{equation}
We will explain at the end of the paper the minor modifications
need to handle the perturbed similarity condition.

 \medskip

 In the remainder of the paper we fix a group of scaling automorphisms 
 $G \sbs \Sim(\cA , \fH , D)$, and let
 $\cA_G = \cA \rtimes G$. We shall also denote by $G_0$ the
 subgroup of isometries in $G$.

 \begin{prop} \label{sst}
 The formula
\begin{equation} \label{tsst}
\s (a \, U) \, = \,\mu (U)^{-1}\, a \, U  , \quad \fl \, U \in G , 
\, a \in \cA ,
\end{equation} 
defines an automorphism  $\s : \cA_G \ra \cA_G$, and 
 $(\cA_G, \fH , D, \s)$ is a $\s$-spectral triple.
Moreover,
 \begin{equation} \label{comms}
[D, a U]_\s  \, = \, [D, a] \, U .
\end{equation}
 \end{prop}
 
 \proof Indeed, one has 
 for any monomials $\, a \, U , \, \, b \, V \in \cA_G = \cA \ltimes G$,
 \begin{eqnarray*}
&& \s(a \, U)\,\s (b \, V)\,=\,
\big( \mu (U)^{-1}\,  \,a \, U \big) \,\big( \mu (V)^{-1}\,  b \, V \big) \, = \\
&&\quad =\, \mu(U V)^{-1} \, a\, (U \, b\, U^\ast ) \, U V\, =\,
  \s\big(a\, U (b))  \, U V \big) \, = \,
 \s(a \, U \,\, b \, V) \,.
 \end{eqnarray*}
Furthermore,  in view of \eqref{sim},
\begin{eqnarray*} 
[D, a \,U]_\s &=& [D, a]\,  U 
\, + \, a \Big( D \, -  \, \mu (U)^{-1}\, U \,D \, U^\ast\Big) \, U \,
\, =\,  [D, a]\,  U \in \cL (\fH) .
\end{eqnarray*}
\endproof 
\medskip

The resulting 
$\s$-spectral triple $(\cA_G, \fH , D, \s)$ will be called {\em twisted
by scaling automorphisms}. 
With the goal of establishing the validity of the Ansatz for twisted spectral triples
of this form, we start with the assumption that the
base spectral triple $(\cA , \fH , D)$ is $(p, \infty)$-summable, and satisfies 
the smoothness condition \eqref{R}. We
recall that the corresponding 
 algebra of pseudodifferential operators  $\Psi^\bullet (\cA  , \fH , D)$ 
 is $\Z$-filtered (by the order) and
$\Z_2$-graded (even/odd), and that it  
includes the subalgebra of  differential operators
$\cD (\cA  , \fH , D)$.

\begin{prop}
With the above notation and hypotheses, 
\begin{enumerate}
\item the action of $\, G$ extends to an action
by automorphisms on $\Psi (\cA  , \fH , D)$, 
\begin{equation} \label{extact}
 P \,\mapsto \, U \trr P := U \, P \, U^\ast , \qquad \fl \, P \in \Psi (\cA  , \fH , D) , \quad U \in G ,
\end{equation} 
which respects both the order filtration and the even/odd grading;
\item  the automorphism $\s \in \Aut ( \cA_G)$ extends to an automorphism 
$\s$ of  the  crossed product algebra $\quad \Psi (\cA_G , \fH , D) := 
  \Psi (\cA  , \fH , D) \rtimes G $, by setting
\begin{equation} \label{tpsi}
\s (P \, U) \, = \,\mu (U)^{-1}\, P \, U  , \quad \fl \, U \in G , 
\, P \in \Psi (\cA  , \fH , D) ,
\end{equation} 
\item the twisted commutators by $D$, $\vert D \vert $ and $D^2$ define twisted
derivations $d_\s$, $\d_\s$, resp. $\nb_\s$,
 of the algebra   $\Psi (\cA_G  , \fH , D) $.
\end{enumerate}
\end{prop}

\proof  The condition \eqref{sim} ensures that $\cD (\cA  , \fH , D)$ remains invariant
under conjugation by $U \in G$, and also implies that 
\begin{equation} \label{zconj}
 U\, \vert D \vert^z\, U^\ast  \, = \, \mu(U)^z \, \vert D \vert^z , \qquad \fl \, z \in \C .
\end{equation}
The verification of the other claims is straightforward.
\endproof
\medskip

We now add the {\em extended simple dimension spectrum} hypothesis: 
{\em there exists
a discrete set $\Sigma_G  \sbs \C$, such that the holomorphic functions}
\begin{equation} \label{zeta2}
  \zeta_B (z) \, = \,  \Tr (B \, |D|^{-z} ) \, , \qquad  \Re z > p \, , \qquad  \fl \, B \in \cB_G :=
  \cB \rtimes G ,
\end{equation}
{\it admit meromorphic extensions to $\C$ with simple poles in $\Sigma_G$}, and
 {\it the functions   
$\, \G (z) \,  \zeta_B (z)$  decay rapidly on finite vertical strips}.

The proof of Proposition \ref{trinv} applies verbatim and shows that
residue functional
\begin{equation*} 
\cutint \cP \,: = \, {\rm Res}_{z=0} \,  \zeta_{\cP} (2z)\, , \qquad 
\cP \in \Psi (\cA_G , \fH , D)
   \end{equation*}
  is automatically a trace. We require it to be $\s$-invariant:
    \begin{equation}  \label{selberg11}
 \cutint\, \s(\cP) \, = \, \cutint\, \cP   .
   \end{equation} 
 This axiom {\em de facto} enforces the Selberg Principle, since it implies
  \begin{equation}  \label{selberg1}
 \cutint P \, U \, = \, 0 , \qquad
\text{if} \quad \mu(U) \,\neq\,1 , \qquad  P \in \Psi (\cA  , \fH , D) , \, \quad U \in G \, ;
   \end{equation}
in particular, the residue functional is necessarily supported on 
$\Psi (\cA_{G_0} , \fH , D) $, where $G_0 := \Isom (\cA , \fH , D)$.

 \medskip
 
\subsection{Twisted JLO brackets} \label{qJLO}

We define the {\em twisted JLO bracket of order} $q$ as the $q+1$-linear
form on $\Psi (\cA_G , \fH , D) $ which for
$\a_0, \ldots , \a_q \in \Psi (\cA  , \fH , D)$ and $U_0, \ldots , U_q \in G$
has the expression
\begin{eqnarray} \label{Jbra1}
&& \quad  \langle \a_0 U_0^\ast , \ldots, \a_q U_q^\ast \rangle_{D} = \int_{\D_q} \Tr 
  \big(\g \, \a_0U_0^\ast \, e^{- s_0 \mu (U_0)^2 D^2} \, \a_1 U_1^\ast \, 
  e^{- s_1 \mu (U_0 U_1)^2 D^2} \cdots \\ \notag
&& \qquad \qquad \qquad   \qquad \qquad   \qquad \qquad  
\cdots \a_qU_q^\ast \, e^{- s_q \mu (U_0 \cdots U_q)^2 D^2} \big)  , 
  \end{eqnarray}  
  where the integration is over the $q$-simplex
$$
\Delta_q:= \bigsetdef{s=(s_0, \ldots ,s_q)\in \R^{q+1}}{s_j\ge 0, \quad s_0+ \ldots +s_q=1} .
$$
\medskip

Throughout the rest of this subsection,
we shall assume that $\a_0, \ldots , \a_q$ are polynomial expressions
in $D$ and the elements of
$\cA$,  $[D, \cA]$, and are homogeneous in $\lb$ as when
$D$ is replaced by $\lb D$.  Given
a JLO bracket as in \eqref{Jbra1}, for any $\ve > 0$ we denote by
$\,   \langle \a_0 U_0^\ast , \ldots, \a_q U_q^\ast \rangle_{D}(\ve) $ the 
expression obtained by replacing every $D$ occurring in each
$\a_0, \ldots , \a_q$
 by $\ve^{1/2} D$.
Equivalently,
\begin{eqnarray} \label{ebra}
  \langle \a_0 U_0^\ast , \ldots, \a_q U_q^\ast \rangle_{D}(\ve) \,=\,  
  \ve^{\frac{m}{2}}   \langle \a_0 U_0^\ast , \ldots, \a_q U_q^\ast \rangle_{\ve^{1/2} D} \, , 
  \end{eqnarray}  
 where $m$ is the total degree in $\lb$ after replacing every
$D$ by $\lb D$ in the product
 $\a_0 \cdots  \a_q$.
As before, by  $U_0, \ldots, U_q$ we denote arbitrary elements in $G$. 

\begin{prop}  \label{asybra} Let $\a_0 \in \cA$, and $\a_1, \ldots , \a_q \in [D, \cA]$.
There is an asymptotic expansion of the form
\begin{equation} \label{asymp}
 \langle \a_0 U_0^\ast , \ldots, \a_q U_q^\ast \rangle_{D} (\ve) \, \sim_{\ve \searrow 0} \,
   \sum_{j \in J} (c_j + c'_j \log \ve)\, \ve^{\frac{q}{2}  - \rho_j}  \, + \, O(1) ,
 \end{equation} 
with $\rho_{0}, \ldots , \rho_{m}$ a finite set of points in the
half-plane $\Re z \geq \frac{q}{2}$.
  \end{prop}

\proof Moving all the unitaries $U_i$ to the rightmost position,
 \begin{eqnarray*} 
&&  \Tr \Big(\g \, \a_0U_0^\ast \, e^{- s_0 \mu (U_0)^2 D^2} U_0 
(U_0^\ast U_1^\ast) \, \a_1 \, (U_1 U_0) U_0^\ast U_1^\ast 
e^{- s_1 \mu (U_0 U_1)^2 D^2} (U_1 U_0)  \cdots \\ 
 &&(U_0^\ast \cdots U_{q-1}^\ast) \, \a_q \, (U_{q-1} \cdots U_0) (U_0^\ast \cdots U_q^\ast)
  e^{- s_q \mu (U_0 \cdots U_q)^2 D^2} (U_q \cdots U_0) U_0^\ast  \cdots U_q^\ast  \Big) ,
 \end{eqnarray*}
 the twisted JLO bracket relative to $\ve^{1/2} D$ takes the form
 \begin{eqnarray} \label{Jbra2}
&& \langle \a_0 U_0^\ast , \ldots, \a_q U_q^\ast \rangle_{\ve^{1/2} D} \, = \\  \notag
&& \quad =\, \int_{\D_q}
  \Tr \big(\g \, \a_0\, e^{- s_1 \ve D^2} \, \a'_1 \, e^{- (s_2-s_1) \ve D^2} \cdots   
\,\a'_q \,
  e^{-(1- s_q) \ve D^2} \, U_0^\ast  \cdots U_q^\ast \big),
\end{eqnarray}
where $\a'_1 = U_0^\ast U_1^\ast \trr \a_1 , \ldots , \a'_q =U_0^\ast \cdots U_{q-1}^\ast \trr \a_q$.

We next use the expansion 
  \begin{equation} \label{EDF}
e^{- \ve D^2} \a \, \sim_{\ve \searrow 0} \, \sum_{n=0}^{\infty} 
\frac{(-1)^n  \ve^n }{n!} \nb^n (\a)\,  e^{-\ve D^2} \, , \qquad \a \in \cD (\cA  , \fH , D) ,
\end{equation}
which is the heat operator analogue of the expansion \eqref{FR2} (\cf also 
 \eqref{IDF} {\em infra}, for
its twisted version), to move the heat operators in Eq. \eqref{Jbra2} to the right
and bring them all in the last position. One obtains
 \begin{eqnarray} \label{Jbra3}
&& \qquad  \langle \a_0 U_0^\ast , \ldots, \a_q U_q^\ast \rangle_{\ve^{1/2} D} \, 
\sim_{\ve \searrow 0} 
\frac{1}{q!} \sum_{N \geq 0} \, \, \sum_{n_1 +\ldots+ n_q =N} 
 \frac{(-1)^N \ve^N }{n_1! \cdots  n_q!} \\ \notag
&&\cdot  \int_{0 \leq  s_1 \leq \ldots \leq  s_q \leq1} s_1^{n_1} \cdots s_q^{n_q}
\cdot  \Tr \big(\g \a_0  \nb^{n_1} (\a'_1)  \cdots   
 \nb^{n_q}(\a'_q)
  e^{-  \ve D^2} U_0^\ast  \cdots U_q^\ast \big) \\ \notag
  && =\,  \frac{1}{q!} \sum_{N \geq 0} \, \, \sum_{n_1 +\ldots+ n_q =N} 
  \frac{(-1)^N \ve^N }{n_1! \cdots n_q! (n_1+1) \cdots (n_1+ \cdots n_q+q)} \cdot\\ \notag
  && \qquad \qquad \cdot   \Tr \big(\g \a_0  \nb^{n_1} (\a'_1)  \cdots   
 \nb^{n_q}(\a'_q)
  e^{-  \ve D^2} U_0^\ast  \cdots U_q^\ast \big).
\end{eqnarray}
 The {\em extended simple dimension spectrum} hypothesis ensures that the
 zeta functions 
 \begin{equation*} \label{zetaN}
 \z_N (z) \, = \,  \Tr \big(\g \,  U_0^\ast  \cdots U_q^\ast\, \a_0\,   \nb^{n_1} (\a'_1) \, \cdots   
\, \nb^{n_q}(\a'_q) \, |D|^{-2z -2N}  \big) , \qquad \Re z > \frac{p}{2} ,
   \end{equation*} 
have meromorphic continuation with simple poles. One has
\begin{equation*} \label{mellN}
\z_N (z) =  \frac{1}{\G (z+N)} \int_0^\infty t^{z +N-1} 
  \Tr \big(\g  U_0^\ast  \cdots U_q^\ast\, \a_0  \nb^{n_1} (\a'_1)  \cdots   
 \nb^{n_q}(\a'_q) e^{- t D^2} \big) dt .
 \end{equation*} 
Proceeding as in the proof of \cite[Theorem II.1]{C-M2}, 
one establishes by means of the
 inverse Mellin transform the existence of
an asymptotic expansion 
\begin{eqnarray} \label{asyN}
&&\ve^{N+\frac{q}{2}}\, \Tr \big(\g \,  U_0^\ast  \cdots U_q^\ast\, \a_0\,   \nb^{n_1} (\a'_1) \, \cdots   
\, \nb^{n_q}(\a'_q) \, e^{- \ve D^2} \big) \, \sim_{\ve \searrow 0}\\ \notag
 && \qquad \qquad \qquad \qquad  \,
   \sum_j (c_{N, j} + c'_{N, j}  \log \ve)\, \ve^{\frac{q}{2}-\rho_{N, j}}\, +  \, O(1)\,  ,
 \end{eqnarray} 
where the exponents  $\rho_{N,j}$ are the real parts of the poles of $\z_N (z)$
 in the half-plane
 $\Re z > \frac{q}{2}$.  \endproof

\begin{defn} We define the {\em constant term}
$\,  \langle \a_0 U_0^\ast , \ldots, \a_q U_q^\ast \rangle_{D} \vert_{0} $ as
 the {\em finite part} ${\rm Pf}_0$ in the asymptotic expansion
 $ \langle \a_0 U_0^\ast , \ldots, \a_q U_q^\ast \rangle_{D} (\ve) $; it is given by
 the coefficient $c_0$ when $\displaystyle \rho_0 = \frac{q}{2}$,
 and is $0$ otherwise.
 \end{defn}
   
  \begin{prop}  \label{consterm} 
The constant term  $\quad  \langle \a_0 U_0^\ast , \ldots, \a_q U_q^\ast \rangle_{D} \vert_{0} \,$
satisfies
 \begin{eqnarray}  \label{selberg2}
  \langle \a_0 U_0^\ast , \ldots, \a_q U_q^\ast \rangle_{D} \vert_{0} &=& 0 ,
  \qquad \text{unless} \qquad \mu(U_0 \cdots U_q) \,=\,1 ;  \\ \label{selberg3}
  \langle \s(\a_0 U_0^\ast) , \ldots, \s(\a_q U_q^\ast) \rangle_D \vert_{0}  &=&
    \langle \a_0 U_0^\ast , \ldots, \a_q U_q^\ast \rangle_{D} \vert_{0} .
  \end{eqnarray}  
  \end{prop}
\proof  Up to a numerical factor,  $\, \langle \a_0 U_0^\ast , \ldots, \a_q U_q^\ast \rangle_{D} \vert_{0} $
coincides with the residue $\,  \Res_{z=0}\, \z_N $. 
In view of  the axiom \eqref{selberg1},
$$  \text{if} \qquad \mu(U_0 \cdots U_q) \,\neq\,1\qquad \text{then} \qquad \Res_{z=0}\, \z_N =0 ,
\qquad \fl \,N \geq 0.
$$
 This proves the property \eqref{selberg2}, which in turn readily implies
\eqref{selberg3}.  
\endproof

\bigskip
 
 In order to compute the constant term,
  we shall employ the elementary Duhamel-type commutator formula
\begin{eqnarray}  \label{DF}
\qquad e^{-(\b-\a)\lb^2 D^2} A  &-& A \, e^{-(\b-\a)\lb^2 \mu^2 D^2} \, = \\ \notag
&-& \int_\a^\b  e^{-(s-\a)\lb^2D^2} \lb^2 (D^2 A - \mu^2 A D^2 ) e^{-(\b-s)\lb^2\mu^2D^2} \, ds ,
 \end{eqnarray} 
where  $A \in \cA_G$, $\lb , \mu > 0$  and $ [\a , \b] \sbs \R$. It is obtained by
integrating the identity
$$
\frac{d}{ds} \left(e^{-(s-\a)D^2} \, A \, e^{-(\b-s)\mu^2D^2}  \right)\, = \,
- e^{-(s-\a)D^2} (D^2 A - \mu^2 A D^2 ) e^{-(\b-s)\mu^2D^2} 
$$
and then replacing $D$ by $\lb D$.

By iterating \eqref{DF}, and using the abbreviation  $\, \nb_\mu (A) = D^2 A - \mu^2 A D^2 $,
 one obtains for any $N \in \N$,
  \begin{equation} \label{IDF}
e^{-t^2 D^2} A \, = \, \sum_{k=0}^{N-1} 
\frac{(-1)^k t^{2k}}{k!} \nb_\mu^k (A)\,  e^{-t^2 \mu^2 D^2} \ +
\, R_N (D, A, \mu, t) .
\end{equation}
The remainder is given by the formula
  \begin{eqnarray} \notag
R_N (D, A, \mu, t)&=& (-1)^{N} t^{2N} 
\int_{\D_N}  
e^{-s_1 t^2D^2}  \nb^N_ \mu(A) e^{-(1-s_1)t^2\mu^2D^2} ds_1\cdots ds_N
 \\ \label{RIDF}
&=&\frac{(-1)^{N}  t^{2N}}{(N-1)!} 
\int_0^1 (1-s)^{N-1}
e^{-s t^2D^2} \nb^N_ \mu(A) e^{-(1-s)t^2\mu^2D^2} ds ,
  \end{eqnarray}
 Using the finite-summability assumption, it is easy to estimate
 the above expression and thus show that 
Eq. \eqref{IDF} does provide an asymptotic expansion as $\, t \searrow 0$.

Applying this expansion for a twisted bracket as in Proposition \ref{asybra}
one obtains
 \begin{eqnarray*} 
&&\langle \a_0 U_0^\ast , \ldots, \a_q U_q^\ast \rangle_{tD} = \int_{\Delta_q} \Tr 
  \big(\g \, \a_0U_0^\ast \, e^{- s_1 \mu (U_0)^2 t^2D^2}  \a_1U_1^\ast
  e^{- (s_2-s_1) \mu (U_0 U_1)^2 t^2 D^2}  \\
  &&  \qquad \qquad \qquad \qquad  \cdots   e^{- (s_q-s_{q-1}) \mu (U_0 \cdots U_{q-1})^2 t^2D^2}
  \a_qU_q^\ast \, e^{-(1- s_q) \mu (U_0 \cdots U_q)^2 t^2 D^2} \big)  \\
  &&\sim_{t \searrow 0}  \sum_{N \geq 0} (-1)^N t^{2N}\sum_{n_1 +\ldots+ n_q =N}  \,
  \frac{\mu(U_0)^{2(n_1 + \ldots + n_{q})} \cdots \mu(U_{q-1})^{2n_q} }{n_1! \cdots n_q!}
  \cdot \\
   && \Tr  \big(\g \, \a_0U_0^\ast \, \nb_\s^{n_1} ( \a_1U_1^\ast ) \cdots
   \nb_\s^{n_q} ( \a_qU_q^\ast ) 
    \, e^{- \mu (U_0 \cdots U_q)^2  t^2D^2} \big) 
     \int_{0 \leq  s_1 \leq \ldots \leq  s_q \leq1} s_1^{n_1} \cdots s_q^{n_q} \\
&&=\, \sum_{N \geq 0} (-1)^N t^{2N} \sum_{n_1 +\ldots+ n_q =N} 
  \frac{\mu(U_0)^{2(n_1 + \ldots + n_{q})} \cdots \mu(U_{q-1})^{2n_q} }{n_1! \cdots n_q! (n_1+1)
   \cdots (n_1+ \cdots n_q+q)} \cdot \\
   &&  \qquad \qquad  \qquad \qquad  \qquad \qquad
    \Tr  \big(\g \, \a_0U_0^\ast \, \nb_\s^{n_1} ( \a_1U_1^\ast ) \cdots
   \nb_\s^{n_q} ( \a_qU_q^\ast )  \, e^{- t^2 D^2} \big) .
   \end{eqnarray*}  
    In view of   \eqref{selberg2}, we may assume 
     $\quad \mu (U_0 \cdots U_q)^{-2(n_1 + \ldots + n_{q})} = 1$; 
   multiplying by $\, \mu (U_0 \cdots U_q)^{-2(n_1 + \ldots + n_{q})}$ , we can continue by
    \begin{eqnarray*} \label{heatcollect2}
   &&=   \sum_{N \geq 0} (-1)^N t^{2N} \sum_{n_1 +\ldots+ n_q =N} 
  \frac{  \mu(U_1)^{-2n_1}  \cdots \mu(U_q)^{-2(n_1 + \ldots + n_q)}}{n_1! \cdots n_q! (n_1+1)
   \cdots (n_1+ \cdots n_q+q)} \cdot \\
   &&  \qquad \qquad  \qquad \qquad  \qquad \qquad
    \Tr  \big(\g \, \a_0U_0^\ast \, \nb_\s^{n_1} ( \a_1U_1^\ast ) \cdots
   \nb_\s^{n_q} ( \a_qU_q^\ast )  \, e^{- t^2D^2} \big) \\
   && =\,  \sum_{N \geq 0} (-1)^N t^{2N} \sum_{n_1 +\ldots+ n_q =N} 
  \frac{1}{n_1! \cdots n_q! (n_1+1) \cdots (n_1+ \cdots n_q+q)}  \cdot \\
   &&  \qquad \qquad
    \Tr  \big(\g \,\a_0U_0^\ast \, 
    \nb_\s^{n_1} \big( \s^{-2n_1}(\a_1U_1^\ast) \big) \cdots
   \nb_\s^{n_q}\big(\s^{-2(n_1 + \ldots + n_q)} ( \a_qU_q^\ast )\big)  \, e^{-t^2 D^2} \big) .
  \end{eqnarray*}  
  
  Comparing with the expansion obtained in the proof of Proposition \ref{asybra}, and
  converting the result into a residue  
  via the Mellin transform, we arrive at the following conclusion.
  
  \begin{prop}  \label{cterm}  
 The constant term has the expression
 \begin{eqnarray*} 
 &&\langle \a_0 U_0^\ast , \ldots, \a_q U_q^\ast \rangle_{D} \vert_{0} 
   = \frac{1}{q!} \sum_{n_1,\ldots , n_q \geq 0} 
  \frac{(-1)^{n_1 + \ldots + n_q}  \G \big(n_1 + \ldots + n_q + \frac{q}{2}\big)}
  {n_1! \cdots n_q! (n_1+1) \cdots (n_1+ \cdots n_q+q)} 
   \cdot  \\ \notag
   && 
      \cutint \g \,\a_0U_0^\ast  
    \nb_\s^{n_1} \big( \s^{-2n_1}(\a_1U_1^\ast) \big) \cdots
   \nb_\s^{n_q}\big(\s^{-2(n_1 + \ldots + n_q)} ( \a_qU_q^\ast )\big) \, 
   |D|^{-2(n_1 + \ldots + n_q) -q} .
  \end{eqnarray*}  
  \end{prop}
 
\bigskip

 \subsection{The cocycle identity} To compute coboundaries of the twisted
 brackets in the cyclic bicomplex, one needs
to establish identities similar to those satisfied by
 the usual JLO brackets, \cf \cite[Lemma 2.2]{GetzSzen}.
 
\begin{lem} \label{forB} With the same notation as in the preceding subsection,
one has
 \begin{eqnarray*}
 \langle \a_0 U_0^\ast , \ldots , \a_q U_q^\ast \rangle_{D} \, =
 \, \sum_{k=0}^{q} 
 \langle \a_0 U_0^\ast , \ldots , 1 , \a_k U_k^\ast , 
\ldots ,  a_{q} U_{q}^\ast \rangle_{D} .
\end{eqnarray*}
\end{lem} 

\proof Proceeding as in \cite[loc. cit.]{GetzSzen}, one writes
 \begin{eqnarray*}
 \langle \a_0 U_0^\ast ,&\ldots& , \a_q U_q^\ast \rangle_{D}\,=\, 
\int_0^1 \langle \a_0 U_0^\ast ,\ldots , \a_q U_q^\ast \rangle_{D} \, ds  \,=\, \\
&=& \int_0^1 ds
 \int_{0\leq s_1 \leq \ldots \leq s_q \leq 1} \Tr 
  \big(\g \, \a_0U_0^\ast \, e^{- s_1 \mu (U_0)^2 D^2}  \a_1 U_1^\ast  \cdot \\
&\,& \qquad \qquad  \qquad \cdot \, e^{- (s_2-s_1) \mu (U_0 U_1)^2 D^2}  \cdots 
 \a_qU_q^\ast \, e^{- (1-s_q) \mu (U_0 \cdots U_q)^2 D^2} \big) \,=\,  \\
&=&  \int_{0\leq s \leq s_1 \leq \ldots \leq s_q \leq 1} \Tr 
  \big(\g \, \a_0U_0^\ast \, e^{- s  \mu (U_0)^2 D^2} \cdot 1 \cdot e^{- (s_1-s) \mu (U_0)^2 D^2} 
  \a_1 U_1^\ast  \cdot   \\
&\,& \qquad \qquad  \qquad \cdot \,e^{- (s_2-s) \mu (U_0 U_1)^2 D^2} 
   \cdots \a_qU_q^\ast \, e^{- (1-s_q) \mu (U_0 \cdots U_q)^2 D^2} \big) \, + \\
 &+&  \int_{0\leq s_1 \leq s \leq \ldots \leq s_q \leq 1} \Tr 
  \big(\g \, \a_0U_0^\ast \, e^{- s_1 \mu (U_0)^2 D^2} \a_1 U_1^\ast  \,
  e^{- (s-s_1) \mu (U_0 U_1)^2 D^2}  \\
&\cdot& 1 \cdot  e^{- (s_2-s) \mu (U_0 U_1)^2 D^2} 
   \cdots \a_qU_q^\ast \, e^{- (1-s_q) \mu (U_0 \cdots U_q)^2 D^2} \big) \quad+ \quad \ldots  \\
 &&\ldots \quad + \quad \int_{0\leq s_1 \leq \ldots \leq s_q \leq s \leq1} \Tr
    \big(\g \, \a_0U_0^\ast \, e^{- s_1 \mu (U_0)^2 D^2} \cdots  
    \a_qU_q^\ast \cdot \\
&\,& \qquad \qquad  \qquad \cdot \, e^{- (s-s_q) \mu (U_0 \cdots U_q)^2 D^2} \cdot 
  1 \cdot e^{- (1-s_q) \mu (U_0 \cdots U_q)^2 D^2}  \big) . 
\end{eqnarray*}

 \begin{lem} \label{forb}
 For $j = 1, \ldots , q-1$, one has
  \begin{eqnarray*} 
&& \langle \a_0 U_0^\ast , \ldots,  \a_{j-1} U_{j-1}^\ast  \cdot \a_j U_j^\ast ,
  \a_{j+1} U_{j+1}^\ast , \dots,  \a_q U_q^\ast \rangle_{D} \\
&&\qquad \qquad \qquad
-  \, \langle \a_0 U_0^\ast , \ldots,  \a_{j-1} U_{j-1}^\ast , 
  \a_j U_j^\ast \cdot \a_{j+1} U_{j+1}^\ast , \dots, 
 \a_q U_q^\ast \rangle_{D} \\
  &&\qquad =\,  \langle \s^2(\a_0U_0^\ast),
 \ldots , \s^2(\a_{j-1} U_{j-1}^\ast) , [D^2 , \a_j U_j^\ast]_\s ,
 \a_{j+1} U_{j+1}^\ast ,  \ldots , \a_qU_q^\ast \rangle_D .
 \end{eqnarray*}
 \end{lem}
 
 \proof 
Making use of the commutator formula \eqref{DF}, one can write
 \begin{eqnarray*} 
&& \langle \a_0 U_0^\ast , \ldots,  \a_{j-1} U_{j-1}^\ast  \cdot \a_j U_j^\ast ,
  \a_{j+1} U_{j+1}^\ast , \dots, 
 \a_q U_q^\ast \rangle_{D} \, - 
  \langle \a_0 U_0^\ast , \ldots,  \a_{j-1} U_{j-1}^\ast , \\
 && \qquad \qquad \qquad \qquad \qquad \qquad \qquad  \qquad \qquad 
\qquad \qquad  \a_j U_j^\ast \cdot \a_{j+1} U_{j+1}^\ast , \dots, 
 \a_q U_q^\ast \rangle_{D} \\
  && =  \int_{\Delta_q} \Tr  \Big(\g \a_0U_0^\ast \, e^{- s_1 \mu (U_0)^2 D^2} 
 \cdots \a_{j-1} U_{j-1}^\ast \big( \int_{t_{j-1}}^{t_{j+1}} 
e^{- (s_j-s_{j-1}) \mu (U_0 \cdots U_{j-1})^2 D^2} \cdot \\
&& \qquad \qquad  \mu (U_0 \cdots U_{j-1})^2 [D^2 , \a_j U_j^\ast]_\s
 e^{- (s_{j+1} - s_j) \mu (U_0 \cdots U_{j})^2 D^2} ds_j \big)\cdot \\
&&\qquad \qquad \qquad \qquad \qquad \qquad  \a_{j+1} U_{j+1}^\ast  \cdots 
 \a_qU_q^\ast \, e^{- (1-s_{q}) \mu (U_0 \cdots U_q)^2 D^2} \Big) \\ 
   && =  \int_{\Delta_q} \Tr  \Big(\g \s^2(\a_0U_0^\ast) \, e^{- s_1 \mu (U_0)^2 D^2} 
 \cdots \s^2(\a_{j-1} U_{j-1}^\ast)  \\
&&\big( \int_{t_{j-1}}^{t_{j+1}} 
e^{- (s_j-s_{j-1}) \mu (U_0 \cdots U_{j-1})^2 D^2} [D^2 , \a_j U_j^\ast]_\s
 e^{- (s_{j+1} - s_j) \mu (U_0 \cdots U_{j})^2 D^2} ds_j \big)\cdot \\
&&\qquad \qquad \qquad  \a_{j+1} U_{j+1}^\ast  \cdots
 \a_qU_q^\ast \, e^{- (1-s_{q}) \mu (U_0 \cdots U_q)^2 D^2} \Big) .
 \end{eqnarray*}

 In contrast with the untwisted case, the cyclic symmetry property  
 is no longer exactly satisfied. It only subsists in a weaker form.

\begin{lem}   \label{cyc}
With $m$ denoting the degree in $D$ of the product
$\a_0 \cdots \a_q$, one has 
 \begin{equation} \label{cyc0}
 \langle \a_0 U_0^\ast , \ldots , \a_q U_q^\ast \rangle_{D}  \vert_{0} \, = \,
\langle  \a_1 U_1^\ast , \ldots, \a_q U_q^\ast ,   \s^{-m} (\a_0U_0^\ast)\rangle_{D} \vert_{0} .
\end{equation}
Moreover, if  $ \mu(U_0 \cdots U_q) = 1$, then
 \begin{equation} \label{cyc1}
  \langle \a_0 U_0^\ast , \ldots , \a_q U_q^\ast \rangle_{D} (\ve) \, = \,
\langle  \a_1 U_1^\ast , \ldots, \a_q U_q^\ast ,  
 \s^{-m} (\a_0U_0^\ast)\rangle_{D} (\mu(U_0)^2 \ve) .
 \end{equation}
\end{lem} 

\proof Indeed,  
 \begin{eqnarray*}
 && \langle \a_0 U_0^\ast , \ldots , \a_q U_q^\ast \rangle_{D} ( \mu(U_0)^{-2}\ve)\, = \, \\
&&=\mu(U_0)^{-m}   \ve^{\frac{m}{2}}  \int_{\D_{q}} \Tr 
  \big(\g  \a_0U_0^\ast e^{- s_{q+1} \ve D^2} \a_1 U_1^\ast e^{- s_1 \mu (U_1)^2 \ve D^2} 
   \cdots \\
   && \cdots \a_qU_q^\ast \, e^{- s_q \mu (U_1 \cdots U_q)^2 \ve D^2} \big) \,
  = \, \ve^{\frac{m}{2}}  \int_{\D_{q}} \Tr 
  \big(\g  \a_1 U_1^\ast e^{- s_1 \mu (U_1)^2 \ve D^2} \cdots \\
  && \cdots \a_qU_q^\ast \, e^{- s_q \mu (U_1 \cdots U_q)^2 \ve D^2} 
  \s^{-m} (\a_0U_0^\ast) e^{- s_{q+1} \ve D^2}\big)  ,
  \end{eqnarray*} 
 which under the assumption  $ \mu(U_0 \cdots U_q) = 1$ equals
    \begin{eqnarray*}
  =&&\ve^{\frac{m}{2}}  \int_{\D_{q}} \Tr 
  \big(\g  \a_1 U_1^\ast e^{- s_1 \mu (U_1)^2 \ve D^2} 
   \cdots \a_qU_q^\ast \, e^{- s_q \mu (U_1 \cdots U_q)^2\ve D^2} \cdot \\
&&\cdot  \s^{-m} (\a_0U_0^\ast) \, e^{- s_{q+1}  \mu  (U_0 \cdots U_q)^2\ve D^2} \big) 
\, = \, 
  \langle \a_1 U_1^\ast , \ldots, \a_q U_q^\ast ,   \s^{-m} (\a_0U_0^\ast)\rangle_{D} ( \ve ) .
  \end{eqnarray*}
This proves \eqref{cyc1}, and also implies the equality of the their constant terms.
On the other hand, if $ \mu(U_0 \cdots U_q) \neq 1$, then
both sides of \eqref{cyc0} vanish, \cf \eqref{selberg2}. 
  \endproof
 
\bigskip

We now introduce the twisted version of the JLO cocycles by defining, for
any $q+1$ elements $A_0, \ldots , A_q \in \cA_G $,
 \begin{equation} \label{Jco1} 
  J^q (D) (A_0,\ldots,A_q) \, = \,
  \langle A_0, [D, \s^{-1} (A_1)]_\s , \ldots,  [D, \s^{-q}(A_q)]_\s \rangle_{D} .
\end{equation}
The collection $\{J^q (D) \}_{q = 0, 2, 4, \ldots}$, resp. $\{J^q (D) \}_{q = 1, 3, 5, \ldots}$,
is a cochain in the entire cyclic cohomology bicomplex of $\cA_G$ but,
because of the failure of cyclic symmetry pointed out above,  this cochain is not
 a cocycle. Instead,
we form the $1$-parameter family
\begin{equation} \label{eJLO}
 J^q (\ve^{1/2} D) (A_0,\ldots,A_q):= \ve^{\frac{q}{2}}   
  \langle A_0, [D, \s^{-1} (A_1)]_\s , \ldots,  [D, \s^{-q}(A_q)]_\s \rangle_{\ve^{1/2} D} \,  , 
  \end{equation}  
where $  \, \ve \in \R^+$,
and passing to the constant term we define 
\begin{equation} \label{cJLO}
  \gimel^q (D) (A_0,\ldots,A_q) \,:=\,     
  \langle A_0, [D, \s^{-1} (A_1)]_\s , \ldots ,  [D, \s^{-q}(A_q)]_\s \rangle_ D \vert_0 \, .
  \end{equation}  
According to Proposition \ref{cterm}, it has the explicit form predicted by the Ansatz
    \begin{eqnarray} \label{Jindex} 
\gimel^q (D) (A_0 ,\ldots ,A_q)&=& \sum_{\bf k}  c_{q, {\bf k} } \, 
 \cutint \, \g \, A_0\,  [D , \s^{-2k_1 - 1} (A_1)]_\s^{(k_1)} \ldots   \\ \nonumber
& &\, \ldots 
 [D ,  \s^{-2(k_1 + \ldots + k_{q}) -q} (A_q)]_\s^{(k_q)} \, \vert D \vert^{-2\vert {\bf k} \vert-q} ,
\end{eqnarray}  
which in particular implies that
\begin{equation} \label{vanq}
\gimel^q (D) = 0 , \qquad  \text{for any} \quad  q > p \, ;
\end{equation}
also, by Proposition \ref{asybra}   
  \begin{eqnarray} \label{kermu}
\gimel^q (D) (a_0 U_0^\ast, \ldots, a_q U_q^\ast) \, = \, 0 ,\qquad \text{if} \qquad
 \mu(U_0^\ast \cdots U_q^\ast)  \neq 1 \, .
  \end{eqnarray} 
 Thus,  $\{\gimel^q (D) \}_{q = 0, 2, 4, \ldots}$, resp. $\{\gimel^q (D) \}_{q = 1, 3, 5, \ldots}$, 
 defines a cochain in the $(b, B)$-bicomplex of $\cA_G$, which is 
  supported on the conjugacy classes from $G_0$.

\begin{thm} \label{constcocy} The cochain
 $\, \gimel^\bullet (D) \, $ satisfies the cocycle identity
 \begin{equation} \label{cocyid}
 b \gimel^{q-1} (D)( a_0 U_0^\ast , \ldots , a_q U_q^\ast) \, + \,
  B \gimel^{q+1} (D)( a_0 U_0^\ast , \ldots , a_q U_q^\ast) \, = \, 0 \, .
  \end{equation} 
\end{thm}

\proof  The first stage of the proof will consist in
establishing the identity
 \begin{eqnarray}  \label{b} 
&&   b \gimel^{q-1} (D)( a_0 U_0^\ast , \ldots , a_q U_q^\ast) \, = \,
   \sum_{j=1}^q (-1)^{j-1}  \langle \s(a_0 U_0^\ast) , \ldots,  \\ \notag  
 && [D, \s^{-(j-2)} (a_{j-1} U_{j-1}^\ast )]_\s,
   [D^2, \s^{-j} ( a_{j} U_{j}^\ast)]_\s, [D, \s^{-j-1} (a_{j+1} U_{j+1}^\ast)]_\s,
  \ldots \\ \notag
  && \qquad \qquad  \qquad \qquad  \qquad \qquad  \qquad \qquad  \qquad \qquad
  \ldots, [D , \s^{-q} (a_q U_q^\ast)]_\s \rangle_{D} \vert_0 .
  \end{eqnarray} 
To this end, we compute
 \begin{eqnarray*} 
 && b J^{q-1} (D)( a_0 U_0^\ast , \ldots , a_q U_q^\ast) \, = \, 
  \langle a_0 U_0^\ast \cdot a_1 U_1^\ast , \ldots ,
[D , \s^{-(q-1)} (a_q U_q^\ast)]_\s \rangle_{D} \, + \\
 && \qquad \qquad \qquad+ \, \sum_{j=1}^{q-1}(-1)^{j} 
 \langle a_0 U_0^\ast , \ldots,  [D, \s^{-j}(a_j U_j^\ast \cdot a_{j+1} U_{j+1}^\ast)]_\s ,
  \dots \rangle_{D} + \\ 
 &&\qquad \qquad \qquad  \qquad \qquad    + \,(-1)^{q} \langle a_q U_q^\ast \cdot a_0 U_0^\ast , 
  \ldots , [D , \s^{-(q-1)} (a_{q-1} U_{q-1}^\ast)]_\s  \rangle_{D} \\ 
    &&\qquad \qquad =  \langle a_0 U_0^\ast \cdot a_1 U_1^\ast , \ldots ,
[D , \s^{-(q-1)} (a_q U_q^\ast]_\s \rangle_{D} \, + \\
&& \qquad - \, \langle a_0 U_0^\ast \,  , \, 
a_1 U_1^\ast \cdot [D, \s^{-1} (a_{2} U_{2}^\ast)]_\s ,
  \dots,  [D , \s^{-(q-1)} (a_q U_q^\ast)]_\s \rangle_{D} \, + \\  
 && + \sum_{j=2}^{q-1}(-1)^{j-1} 
 \langle a_0 U_0^\ast , \ldots \\
 && \qquad \qquad \qquad \qquad\ldots,  [D, \s^{-(j-1)} (a_{j-1} U_{j-1}^\ast )]_\s \cdot 
 \s^{-(j-1)} ( a_{j} U_{j}^\ast),
  \ldots,    \rangle_{D} \, +\\ 
 &&+ \sum_{j=2}^{q-1}(-1)^{j} 
 \langle a_0 U_0^\ast , \ldots,  \s^{-(j-1)} (a_j U_j^\ast )\cdot [D, \s^{-j} (a_{j+1} U_{j+1}^\ast)]_\s ,
  \ldots,  \rangle_{D} \, + \\  
 &&\qquad  + (-1)^{(q-1)} \langle a_0 U_0^\ast , [D, \s^{-1} (a_{1} U_{1}^\ast)]_\s , \ldots \\
&&\qquad  \qquad  \ldots,  [D, \s^{-(q-1)} (a_{q-1} U_{q-1}^\ast )]_\s \cdot 
 \s^{-(q-1)} ( a_{q} U_{q}^\ast) \rangle_{D} \\ 
 &&\qquad \qquad  \qquad \qquad   + \,(-1)^{q}  \langle a_q U_q^\ast \cdot a_0 U_0^\ast , 
  \ldots , [D , \s^{-(q-1)} (a_{q-1} U_{q-1}^\ast)]_\s  \rangle_{D} ,
  \end{eqnarray*} 
which by Lemma \ref{forb} is equal to
 \begin{eqnarray*}
    &&\qquad \qquad \qquad  \langle \s^2(a_0 U_0^\ast) , [D^2,  a_1 U_1^\ast]_\s , \ldots ,
  [D , \s^{-(q-1)} (a_q U_q^\ast)]_\s \rangle_{D} \, + \\  
 && + \sum_{j=2}^{q-1}(-1)^{j-1} 
 \langle \s^2(a_0 U_0^\ast) , \ldots,  [D, \s^{-(j-3)} (a_{j-1} U_{j-1}^\ast )]_\s, 
[D^2, \s^{-(j-1)} ( a_{j} U_{j}^\ast)]_\s , \\
&& \qquad \qquad \qquad  \qquad \qquad [D, \s^{-j} (a_{j+1} U_{j+1}^\ast)]_\s,
  \ldots,  [D , \s^{-(q-1)} (a_q U_q^\ast)]_\s \rangle_{D} \\  
 &&+ (-1)^{(q-1)} \langle a_0 U_0^\ast , [D, \s^{-1} (a_{1} U_{1}^\ast)]_\s , \ldots \\
&&\qquad \qquad \qquad  \ldots,  [D, \s^{-(q-1)} (a_{q-1} U_{q-1}^\ast )]_\s \cdot 
 \s^{-(q-1)} ( a_{q} U_{q}^\ast) \rangle_{D} \, +\\ 
 && \qquad  \qquad \qquad \qquad \qquad   + \,(-1)^{q}  \langle a_q U_q^\ast \cdot a_0 U_0^\ast ,
  \ldots , [D , \s^{-(q-1)} (a_{q-1} U_{q-1}^\ast)]_\s  \rangle_{D} .
  \end{eqnarray*} 
  At this point we 
 pass to the constant terms and use Eq. \eqref{cyc0} for the last two terms,  to replace
 them by the sum
 \begin{eqnarray*}
 &&  (-1)^{(q-1)} \langle [D, \s^{-1} (a_{1} U_{1}^\ast)]_\s , 
  \ldots \\
  && \ldots , [D, \s^{-(q-1)} (a_{q-1} U_{q-1}^\ast )]_\s \cdot 
 \s^{-(q-1)} ( a_{q} U_{q}^\ast) ,  \s^{-(q-1)} ( a_0 U_0^\ast) \rangle_{D} \vert_0 \\ 
 &&  \qquad \qquad   +\,  \,(-1)^q \langle [D, \s^{-1} (a_{1} U_{1}^\ast)]_\s , 
  \ldots  [D , \s^{-(q-1)} (a_{q-1} U_{q-1}^\ast)]_\s,  \\
  && \qquad \qquad  \qquad \qquad \qquad \s^{-(q-1)}(a_q U_q^\ast)
   \cdot \s^{-(q-1)}(a_0 U_0^\ast) \rangle_{D}  \vert_0 \, .
  \end{eqnarray*}  
In turn, by Lemma \ref{forb} this equals 
   \begin{eqnarray*}
&&=  (-1)^{(q-1)} \langle [D, \s (a_{1} U_{1}^\ast)]_\s , 
  \ldots,  [D, \s^{-(q-3)} (a_{q-1} U_{q-1}^\ast )]_\s , \\
&& \qquad  \qquad  \qquad  \quad
 [D^2, \s^{-(q-1)} ( a_{q} U_{q}^\ast)]_\s ,  \s^{-(q-1)} ( a_0 U_0^\ast) \rangle_{D}\vert_0 .
 \end{eqnarray*} 
 Applying once again Eq. \eqref{cyc0}, 
 and taking into account that the homogeneity degree in $D$ is $q+1$, 
the above expression becomes
 \begin{eqnarray*}
&&= (-1)^{(q-1)} \langle \s^2(a_{0} U_{0}^\ast) , [D, \s (a_{1} U_{1}^\ast)]_\s , 
  \ldots,  [D, \s^{-(q-3)} (a_{q-1} U_{q-1}^\ast )]_\s , \\
&&\qquad  \qquad  \qquad  [D^2, \s^{-(q-1)} ( a_{q} U_{q}^\ast)]_\s \rangle_{D} \vert_0 .
  \end{eqnarray*}  
 Summing up, one obtains
   \begin{eqnarray*} 
 && b \gimel^{q-1} (D)( a_0 U_0^\ast , \ldots , a_q U_q^\ast) \, = \\
   &&\sum_{j=1}^{q}(-1)^{j-1} 
 \langle \s^2(a_0 U_0^\ast) , \ldots,  [D, \s^{-(j-3)} (a_{j-1} U_{j-1}^\ast )]_\s, 
[D , [D, \s^{-(j-1)} ( a_{j} U_{j}^\ast)]_\s]_\s , \\
&& \qquad \qquad \qquad \qquad \quad  [D, \s^{-j} (a_{j+1} U_{j+1}^\ast)]_\s,
  \ldots,  [D , \s^{-(q-1)} (a_q U_q^\ast)]_\s \rangle_{D} \vert_0 .
  \end{eqnarray*} 
  Combined with the $\s$-invariance property \eqref{selberg3}, this completes the
proof of Eq. \eqref{b}.
\medskip 

The second stage of the proof will show that the $B$-boundary 
  \begin{eqnarray*} 
 && B \gimel^{q+1} (D)( a_0 U_0^\ast , \ldots , a_q U_q^\ast) \, = \\
   && =\, \sum_{k=0}^{q} (-1)^{kq}\, \gimel^{q+1} (D)(1 ,  a_{k}U_{k}^\ast, 
\ldots  a_{q}U_{q}^\ast , a_{0} U_{0}^\ast, \ldots , a_{k-1} U_{k-1}^\ast ) ,
\end{eqnarray*}
 satisfies the identity
 \begin{eqnarray} \label{B}
 B  \gimel^{q+1} (D)( a_0 U_0^\ast , \ldots , a_q U_q^\ast)  &=&
  \langle  [D, \s^{-1}(a_{0} U_{0}^\ast)]_\s , \ldots \\ \notag
 &&\qquad 
  \ldots , [D, \s^{-(q+1)}(a_{q}U_{q}^\ast)]_\s \rangle_{D} \vert_0.
\end{eqnarray}
 
To this end, we note that by Lemma \ref{forB},
 \begin{eqnarray*}
&&\langle  [D, \s^{-1}(a_{0} U_{0}^\ast)]_\s , \ldots , 
 [D, \s^{-(q+1)}(a_{q}U_{q}^\ast)]_\s \rangle_{D} \, =\\
 &&=  \sum_{k=0}^{q} 
 \langle [D, \s^{-1}(a_{0} U_{0}^\ast)]_\s  , \ldots , 1 ,  [D, \s^{-(k+1)}(a_{k}U_{k}^\ast)]_\s, 
\ldots ,   [D, \s^{-(q+1)}(a_{q}U_{q}^\ast)]_\s \rangle_{D} .
\end{eqnarray*}
Passing to the constant term, we apply Eq. \eqref{cyc0} $k$-times to
the $k$-th term of the sum and rewrite it in the form   
 \begin{eqnarray*}
 &&  \langle [D, \s^{-1}(a_{0} U_{0}^\ast)]_\s  , \ldots , 1 ,  [D, \s^{-(k+1)}(a_{k}U_{k}^\ast)]_\s, 
\ldots ,   [D, \s^{-(q+1)}(a_{q}U_{q}^\ast)]_\s \rangle_{D} \vert_0 \\
    && = (-1)^{kq} \langle 1 ,  [D, \s^{-(k+1)}(a_{k}U_{k}^\ast)]_\s, 
\ldots  [D, \s^{-(q+1)}(a_{q}U_{q}^\ast)]_\s , [D, \s^{-(q+2)}(a_{0} U_{0}^\ast)]_\s ,\\
&& \qquad \qquad  \qquad \qquad 
  \ldots , [D, \s^{-(q+k+1)}(a_{k-1} U_{k-1}^\ast)]_\s\rangle_{D}\vert_0 .
\end{eqnarray*}
Summing up, one obtains 
 \begin{eqnarray*}
&& \langle  [D, \s^{-1}(a_{0} U_{0}^\ast)]_\s , \ldots ,  
[D, \s^{-(q+1)}(a_{q}U_{q}^\ast)]_\s \rangle_{D}\vert_0 \\
 && = \, \sum_{k=0}^{q} (-1)^{kq} \langle 1 ,  [D, \s^{-(k+1)}(a_{k}U_{k}^\ast)]_\s, 
\ldots  [D, \s^{-(q+1)}(a_{q}U_{q}^\ast)]_\s , [D, \s^{-(q+2)}(a_{0} U_{0}^\ast)]_\s ,\\
&& \qquad \qquad  \qquad \qquad 
  \ldots , [D, \s^{-(q+k+1)}(a_{k-1} U_{k-1}^\ast)]_\s\rangle_{D}\vert_0  \\
   && = \, \sum_{k=0}^{q} (-1)^{kq}\, \gimel^{q+1} (D)(1 ,  a_{k}U_{k}^\ast, 
\ldots  a_{q}U_{q}^\ast , a_{0} U_{0}^\ast, \ldots , a_{k-1} U_{k-1}^\ast ) ,
\end{eqnarray*}
which proves Eq. \eqref{B}.
 
To relate the two identities satisfied by the coboundary operators, we use
 the generalized Leibniz rule \eqref{tdergen} and write
 \begin{eqnarray*} 
&&   [D, a_0U_0^\ast \, e^{- s_0 \mu (U_0)^2 D^2}  \cdots
[D,  \s^{-q}(a_qU_q^\ast)]_\s \, e^{- s_q \mu (U_0 \cdots U_q)^2 D^2}]_\s \\
&&= [D,  a_0U_0^\ast ]_\s\, e^{- s_0 \mu (U_0)^2 D^2}   \cdots
[D,  \s^{-q}(a_qU_q^\ast)]_\s \, e^{- s_q \mu (U_0 \cdots U_q)^2 D^2} + \\
&& \s( a_0U_0^\ast) \, e^{- s_0 \mu (U_0)^2 D^2} [D^2, \s^{-1}(a_1U_1^\ast)]_\s \,
e^{- s_1 \mu (U_0 U_1)^2 D^2} \cdots  \\
&& \qquad \qquad\qquad \qquad  \cdots
[D,  \s^{-q}(a_qU_q^\ast)]_\s \, e^{- s_q \mu (U_0 \cdots U_q)^2 D^2}\,  + \ldots  +\\
&&(-1)^{q-1} \s(a_0U_0^\ast )e^{- s_0 \mu (U_0)^2 D^2} [D, a_1U_1^\ast]_\s
e^{- s_1 \mu (U_0 U_1)^2 D^2} \cdots \\
&& \qquad \qquad\qquad \qquad \qquad \qquad\qquad \qquad \cdots
[D^2,  \s^{-q}(a_qU_q^\ast)]_\s  e^{- s_q \mu (U_0 \cdots U_q)^2 D^2}.
  \end{eqnarray*}
Since, in view of the Selberg property \eqref{selberg1},  $\displaystyle \cutint $ vanishes on 
twisted graded commutators, one obtains
 \begin{eqnarray*} 
&&- \langle [D, a_0U_0^\ast]_\s, [D,\s^{-1}(a_1U_1^\ast)]_\s , \cdots ,
[D,  \s^{-q}(a_qU_q^\ast)]_\s \rangle_D \vert_0\, = \\
  &&= \sum_{j=1}^q (-1)^{j-1} 
 \langle \s(a_0 U_0^\ast) , \ldots,  [D, \s^{-(j-2)} (a_{j-1} U_{j-1}^\ast )]_\s, 
[D^2, \s^{-j} ( a_{j} U_{j}^\ast)]_\s , \\
&&\qquad \qquad \qquad \qquad \qquad \qquad \quad  [D, \s^{-j-1} (a_{j+1} U_{j+1}^\ast)]_\s,
  \ldots,  [D , \s^{-q} (a_q U_q^\ast)]_\s \rangle_{D} \vert_0 .
   \end{eqnarray*} 
We can now rewrite Eq. \eqref{b} in the form
   \begin{eqnarray*} 
 && b \gimel^{q-1} (D)( a_0 U_0^\ast , \ldots , a_q U_q^\ast) \, = \\
   &&\qquad \qquad -\, \langle [D, a_0U_0^\ast]_\s, [D,\s^{-1}(a_1U_1^\ast)]_\s , \cdots ,
[D,  \s^{-q}(a_qU_q^\ast)]_\s \rangle_D \vert_0 \, ,
  \end{eqnarray*} 
  Using once more the invariance property
 \eqref{selberg3} 
and comparing with Eq. \eqref{B} one obtains the desired cocycle identity.
\endproof

 \medskip
 
\subsection{Transgression and proof of the Ansatz }  \label{TqJLO}
In  view of the the property \eqref{kermu},  we can restrict
our considerations  
to the $(b, B)$-subcomplex $CC_{G_0}^\ast (\cA_G)$ of
cochains supported by the conjugacy classes in $G_0$. Far from being a
mere convenience, this restriction is actually essential for the validity of
the ensuing calculations.

 We denote by $\i (V)$ the {\it twisted contraction} operator on $CC^\ast (\Psi_G)$,
 \begin{eqnarray*}   
   \i_\s  (V)  \langle A_0 , &\ldots&, A_q \rangle_{D} \, = \\
 &=&  \sum_{k=0}^q (-1)^{(\#A_0 + \ldots + \#A_k) \#V}
 \langle \s^2( A_0), \ldots , \s^2( A_k) ,
 V ,  A_{k+1} ,\ldots,  A_q \rangle_{D} , 
\end{eqnarray*}
where $\#A$ stands for the degree of $A$, and extend it to 
cochain-valued functions by setting
 \begin{equation*}   
 \i_\s (V)  \langle A_0 , \ldots, A_q \rangle_{D} (\ve) := 
   \iota (V) \big( \langle  A_0 , \ldots, A_q \rangle_{D} (\ve)\big) , \qquad \ve \in \R^+ .
\end{equation*}

In what follows, we shall denote by $\tau \mapsto D_\tau$
 one of the following two families of operators 
$ D_t = tD , \, t \in \R^+$ and $D_u = D |D|^{-u} , \, u \in [0, 1]$, and will denote by
$\dot{D}$ the corresponding derivative.
In each case,  we define the 
cochains $\Jb^q (D_\tau, V)  \in CC_{G_0}^\ast (\cA_G)$
by the formula
 \begin{eqnarray} \label{Jbco1} 
&&\Jb^q (D_\tau, V) ( a_0 U_0^\ast , \ldots , a_q U_q^\ast) \, = \\  \notag
&& \qquad \qquad = \, 
   \i_\s  (V)  \langle a_0 U_0^\ast , [D_\tau, \s^{-1}(a_1 U_1^\ast)]_\s ,\ldots, 
 [D_\tau, \s^{-q}(a_q U_q^\ast)]_\s  \rangle_{D_\tau} \, ,
\end{eqnarray}
where $V$ will be either the (odd) operator $\dot{D_\tau}$ or the (even)
operator $[D_\tau, \dot{D_\tau}]$.

We would like to evaluate the expression
 \begin{equation} \label{Jbco2} 
\frac{d}{d\tau} J^q (D_\tau)  \, + \,b\Jb^{q-1} (D_\tau, \dot{D_\tau}) + 
  B \Jb^{q+1} (D_\tau, \dot{D_\tau}) \, ,
\end{equation}
which vanishes in the untwisted case (\cf \eg \cite[Prop. 10.12]{Var}). 
The derivative   
 \begin{eqnarray*} 
&& \frac{d}{d\tau} J^q (D_\tau) (a_0 U_0^\ast, \ldots, a_q U_q^\ast)\,=\,  \int_{\D_q}  \frac{d}{d\tau} 
 \Tr \Big(\g a_0 U_0^\ast e^{- s_1 \mu(U_0)^2  D_\tau^2}[D_\tau,  \s^{-1}(a_1 U_1^\ast)]_\s \\
&& \qquad  e^{- (s_2-s_1) \mu(U_0 U_1)^2 D_\tau^2} \cdots [D_\tau, \s^{-q}(a_q U_q^\ast)]_\s 
  e^{- (1- s_q)  \mu(U_0 \cdots U_q)^2 D_\tau^2} \Big) .
\end{eqnarray*}
splits into two sums of terms. The first sum simply consists of the derivatives of
the twisted commutators
 \begin{equation} \label{K1}
 K^q (D_\tau) (a_0 U_0^\ast, \ldots, a_q U_q^\ast)\,:=\,
 \sum_{j=1}^q \langle a_0 U_0^\ast, \ldots , [\dot{D}_\tau, \s^{-j}(a_j U_j^\ast)]_\s ,
\ldots \rangle_{D_\tau} .
\end{equation}
To evaluate the second sum
one relies, as in the standard case, on the Duhamel formula
 \begin{equation*}
 \frac{d}{d\tau} e^{-D_\tau^2}   \, = \, - \,
 \int_0^1  e^{- s D_\tau^2}  [D_\tau, \dot{D}_\tau] \, e^{-(1-s) D_\tau^2} \, ds .
 \end{equation*}
 By applying it in the form
 \begin{equation} \label{DF2}
 \frac{d}{d\tau} e^{-(s_{j+1} - s_j)\mu^2 D_\tau^2}   \, = \, - \, \mu^2
 \int_{s_j}^{s_{j+1}}  e^{- (s-s_j)\mu^2 D_\tau^2}  [D_\tau, \dot{D}_\tau] \, 
 e^{-(s_{j+1}-s) \mu^2D_\tau^2} \, ds \, ,
 \end{equation}
one obtains
 \begin{eqnarray*} 
&&  \sum_{j=0}^q \int_{\D_q}   \Tr \Big(\g  \cdots  [D_\tau, \s^{-j}(a_j U_j^\ast)]_\s 
 \frac{d}{d\tau}  e^{- (s_{j+1}-s_j)  \mu(U_0 \cdots U_j)^2  D_\tau^2} \cdots \Big)= \\
&&   \qquad  = -  \sum_{j=0}^q   \mu(U_0 \cdots U_j)^2 \int_{\D_{q+1}}  
 \Tr \Big(\g  \cdots  [D_\tau, \s^{-j}(a_j U_j^\ast)]_\s \\
&& \qquad  \qquad   e^{- (s-s_j)  \mu(U_0 \cdots U_j)^2 D_\tau}  [D_\tau, \dot{D}_\tau] \, 
   e^{-(s_{j+1}-s)  \mu(U_0 \cdots U_j)^2 D_\tau^2}
 \cdots \Big)  = \\
 &&= -\sum_{j=0}^q  \int_{\D_{q+1}}    
 \langle \s^2 (a_0 U_0^\ast), \ldots ,  [D_\tau, \s^{-(j-2)}(a_j U_j^\ast)]_\s ,
[D_\tau, \dot{D}_\tau] , \ldots \\
&&\ldots , [D_\tau, \s^{-q}(a_q U_q^\ast)]_\s \rangle_{D_\tau} \qquad
  = \qquad -  \Jb^q (D_\tau, [D_\tau, \dot{D_\tau}]) ( a_0 U_0^\ast , \ldots , a_q U_q^\ast) \, .
\end{eqnarray*}
This gives the identity
\begin{equation} \label{derJ}
 \frac{d}{d\tau} J^q (D_\tau) \, = \,  K^q (D_\tau) \, - \,  \Jb^q (D_\tau, [D_\tau, \dot{D_\tau}]) .
\end{equation}

On the other hand, in order to
evaluate the coboundary of $\Jb^\bullet (D_\tau, \dot{D_\tau}) $,  
 as in the proof of Theorem \ref{constcocy}, we apply 
 the Leibniz rule to the integrand for the expression of
$ \Jb^q (D_\tau, \dot{D_\tau}) ( a_0 U_0^\ast , \ldots , a_q U_q^\ast)$.
By abuse of notation, we write this bracket operation in the form
 \begin{equation*} 
[D_\tau,\,  \Jb^q (D_\tau, \dot{D_\tau}) ( a_0 U_0^\ast , \ldots , a_q U_q^\ast)]_\s  \, ,
\end{equation*}
 and compute it as follows
 \begin{eqnarray*}
&& [D_\tau,   \i_\s  (\dot{D_\tau})  \langle a_0 U_0^\ast , [D_\tau, \s^{-1}(a_1 U_1^\ast)]_\s ,\ldots, 
 [D_\tau, \s^{-q}(a_q U_q^\ast)]_\s  \rangle_{D_\tau}]_\s\, = \\
 &&=  [D_\tau,    \langle \s^2(a_0 U_0^\ast) , \dot{D_\tau} ,
 [D_\tau, \s^{-1}(a_1 U_1^\ast)]_\s ,\ldots, 
 [D_\tau, \s^{-q}(a_q U_q^\ast)]_\s  \rangle_{D_\tau}]_\s\, + \ldots \\
  &&=  \langle [D_\tau, \s^2(a_0 U_0^\ast)]_\s ,  \dot{D_\tau} , 
   [D_\tau, \s (a_1 U_1^\ast)]_\s , \ldots, 
 [D_\tau, \s^{-q}(a_q U_q^\ast)]_\s  \rangle_{D_\tau} \\
 && \, +  \langle \s^4(a_0 U_0^\ast) , [D_\tau, \dot{D_\tau}] , 
  [D_\tau, \s^{-1}(a_1 U_1^\ast)]_\s ,\ldots , 
 [D_\tau, \s^{-q}(a_q U_q^\ast)]_\s  \rangle_{D_\tau} + \\
 &&+   \langle \s^4(a_0 U_0^\ast) ,  \dot{D_\tau} , 
   [D_\tau, \s (a_1 U_1^\ast)]_\s ,  [D_\tau, \s^{-2}(a_2 U_2^\ast)]_\s ,\ldots, 
 [D_\tau, \s^{-q}(a_q U_q^\ast)]_\s  \rangle_{D_\tau} \\
 && \qquad  \qquad   \qquad  + \ldots  \text{and so on}.
\end{eqnarray*}
 There are two kinds of terms appearing in this sum.
Those which contain the term
$[D_\tau ,  \s^{2}(a_0U_0^\ast)]_\s$ come from 
 \begin{eqnarray*} 
&& \iota(\dot{D_\tau}) \, \langle [D_\tau, a_0U_0^\ast]_\s, [D_\tau,
\s^{-1}(a_1U_1^\ast)]_\s , \cdots ,
[D_\tau,  \s^{-q}(a_qU_q^\ast)]_\s \rangle_{D_\tau} ,
 \end{eqnarray*}
 and they closely resemble those appearing in 
$B  \Jb^{q+1} (D_\tau, \dot{D_\tau})$. 
The remaining terms are of the form 
  \begin{eqnarray*} 
  && \langle \s^4(a_0 U_0^\ast) , \ldots,  \dot{D}_\tau , \ldots, 
  [D_\tau, \s^{-(j-2)} (a_{j-1} U_{j-1}^\ast )]_\s, 
[D_\tau^2, \s^{-j} ( a_{j} U_{j}^\ast)]_\s, \\
&&\qquad \qquad \qquad \qquad \qquad \qquad \quad  [D_\tau, \s^{-j-1} (a_{j+1} U_{j+1}^\ast)]_\s,
  \ldots,  [D_\tau , \s^{-q} (a_q U_q^\ast)]_\s \rangle_{D_s} ,
   \end{eqnarray*} 
or
  \begin{eqnarray*} 
  && \langle \s^4(a_0 U_0^\ast) , \ldots,  
  [D_\tau, \s^{-(j-2)} (a_{j-1} U_{j-1}^\ast )]_\s, [D_\tau^2, \s^{-j} ( a_{j} U_{j}^\ast)]_\s, \ldots \\
&&\qquad \qquad \ldots,  \dot{D}_\tau , \ldots,  [D_\tau, \s^{-j-1} (a_{j+1} U_{j+1}^\ast)]_\s,
  \ldots,  [D_\tau , \s^{-q} (a_q U_q^\ast)]_\s \rangle_{D_\tau} ,
   \end{eqnarray*} 
and they match those occurring in
$b  \Jb^{q-1} (D_\tau, \dot{D_\tau})+ K^q$, plus terms which contain
$[D_\tau , \dot{D}_\tau ]$ and account for 
$ \Jb^{q-1} (D_\tau, [D_\tau , \dot{D}_\tau ])$.  

Indeed, the $b$-coboundary of $\Jb^\bullet (D_\tau, \dot{D_\tau}) $,
we write
  \begin{eqnarray*} 
 && b \Jb^{q-1} (D_\tau, \dot{D_\tau}) ( a_0 U_0^\ast , \ldots , a_q U_q^\ast) \, =  \\
    &&\qquad = \,  \i_\s (\dot{D_\tau})
 \langle a_0 U_0^\ast \cdot a_1 U_1^\ast ,  [D_\tau, \s^{-1}(a_2 U_2^\ast)]_\s ,
  \ldots,  [D_\tau, \s^{-q+1}(a_q U_q^\ast)]_\s \rangle_{D_\tau}   \\
&& \qquad - \,  \i_\s (\dot{D_\tau}) \langle a_0 U_0^\ast \,  , \, 
a_1 U_1^\ast \cdot [D_\tau, \s^{-1} (a_{2} U_{2}^\ast)]_\s ,
  \dots,  [D_\tau , \s^{-(q-1)} (a_q U_q^\ast)]_\s \rangle_{D_\tau}  \,  + \\  
 && \sum_{j=2}^{q-1}(-1)^{j-1}   \i_\s (\dot{D_\tau})
 \langle a_0 U_0^\ast , \ldots,  [D_\tau, \s^{-(j-1)} (a_{j-1} U_{j-1}^\ast )]_\s \cdot 
 \s^{-(j-1)} ( a_{j} U_{j}^\ast),
  \ldots \rangle_{D_\tau}  \\ 
 &&+ \sum_{j=2}^{q-1}(-1)^{j}   \i_\s (\dot{D_\tau})
 \langle a_0 U_0^\ast , \ldots,  \s^{-(j-1)} (a_j U_j^\ast )\cdot [D_\tau, \s^{-j} (a_{j+1} U_{j+1}^\ast)]_\s ,
  \ldots \rangle_{D_\tau}\,  + \\  
 && + (-1)^{(q-1)}   \i_\s (\dot{D_\tau}) \langle a_0 U_0^\ast , 
 [D_\tau, \s^{-1} (a_{1} U_{1}^\ast)]_\s , \ldots \\
 &&\qquad \qquad \qquad     \ldots,  [D_\tau, \s^{-(q-1)} (a_{q-1} U_{q-1}^\ast )]_\s \cdot 
 \s^{-(q-1)} ( a_{q} U_{q}^\ast) \rangle_{D_\tau} + \\ 
 &&\qquad \qquad \qquad     + \,(-1)^{q}   \i_\s (\dot{D_\tau})
 \langle a_q U_q^\ast \cdot a_0 U_0^\ast , 
  \ldots , [D_\tau , \s^{-(q-1)} (a_{q-1} U_{q-1}^\ast)]_\s  \rangle_{D_\tau}\, . 
  \end{eqnarray*}  
  Let us take a closer look at the first two terms, and expand $ \i_\s (\dot{D_\tau})$. One has   
  \begin{eqnarray*} 
     && \i_\s (\dot{D_\tau})
 \langle a_0 U_0^\ast \cdot a_1 U_1^\ast ,  [D_\tau, \s^{-1}(a_2 U_2^\ast)]_\s ,
  \ldots,  [D_\tau, \s^{-q+1}(a_q U_q^\ast)]_\s \rangle_{D_\tau} \,   \\
&& \qquad - \,   \i_\s (\dot{D_\tau}) \langle a_0 U_0^\ast \,  , \, 
a_1 U_1^\ast \cdot [D, \s^{-1} (a_{2} U_{2}^\ast)]_\s ,
  \dots,  [D_\tau , \s^{-(q-1)} (a_q U_q^\ast)]_\s \rangle_{D_\tau}  \,  \\  
&&=    \langle \s^2(a_0 U_0^\ast \cdot a_1 U_1^\ast) , \, 
\dot{D_\tau} , \, [D_\tau, \s^{-1}(a_2 U_2^\ast)]_\s ,
  \ldots,  [D_\tau, \s^{-q+1}(a_q U_q^\ast)]_\s \rangle_{D_\tau} \,   \\
  && \qquad - \, \langle \s^2(a_0 U_0^\ast ), \, \dot{D_\tau} , \, 
a_1 U_1^\ast \cdot [D_\tau, \s^{-1} (a_{2} U_{2}^\ast)]_\s ,
  \dots,  [D_\tau , \s^{-(q-1)} (a_q U_q^\ast)]_\s \rangle_{D_\tau}  \,  \\  
  &&+ \sum_{k=2}^q  (-1)^k
 \langle \s^2(a_0 U_0^\ast \cdot a_1 U_1^\ast) ,  [D_\tau, \s (a_2 U_2^\ast)]_\s , \ldots , \\
&& [D_\tau, \s^{-(k-3)}(a_k U_k^\ast)]_\s , 
\dot{D_\tau} , [D_\tau, \s^{-k}(a_{k+1} U_{k+1}^\ast)]_\s , 
  \ldots,  [D_\tau, \s^{-(q-1)}(a_q U_q^\ast)]_\s \rangle_{D_\tau} \,  \\
 && - \sum_{k=2}^q   (-1)^k \langle \s^2(a_0 U_0^\ast), 
\s^2(a_1 U_1^\ast) \cdot [D_\tau, \s(a_{2} U_{2}^\ast)]_\s , \ldots ,  \\
&&  [D_\tau, \s^{-(k-3)}(a_k U_k^\ast)]_\s , \dot{D_\tau} , 
[D_\tau, \s^{-k}(a_{k+1} U_{k+1}^\ast)]_\s , 
  \ldots,  [D_\tau, \s^{-(q-1)}(a_q U_q^\ast)]_\s \rangle_{D_\tau} \,  \\
  &&=   \langle \s^2(a_0 U_0^\ast \cdot a_1 U_1^\ast) , \, \dot{D_\tau} , \, 
  [D_\tau, \s^{-1}(a_2 U_2^\ast)]_\s ,
  \ldots,  [D_\tau, \s^{-q+1}(a_q U_q^\ast)]_\s \rangle_{D_\tau}   \\
&& \qquad  - \,  \langle \s^2(a_0 U_0^\ast), \, \dot{D_\tau} , \, 
a_1 U_1^\ast \cdot [D_\tau, \s^{-1} (a_{2} U_{2}^\ast)]_\s ,
  \dots,  [D_\tau , \s^{-(q-1)} (a_q U_q^\ast)]_\s \rangle_{D_\tau} \\  
  &&+ \sum_{k=2}^q (-1)^k
 \langle \s^4(a_0 U_0^\ast), [D_\tau^2, \s^2(a_1 U_1^\ast)]_\s ,  [D_\tau, \s (a_2 U_2^\ast)]_\s ,  \ldots , \\
&&  [D_\tau, \s^{-(k-3)}(a_k U_k^\ast)]_\s , \dot{D_\tau} , 
[D_\tau, \s^{-k}(a_{k+1} U_{k+1}^\ast)]_\s , 
  \ldots,  [D_\tau, \s^{-(q-1)}(a_q U_q^\ast)]_\s \rangle_{D_\tau} \, 
   \end{eqnarray*} 
  where  we have used Lemma \ref{forb}  after the first pair of terms.
  
  We now look at pairs of terms indexed by the same $j=2, \ldots , q-1$,
  \begin{eqnarray*}
  && \i_\s (\dot{D_\tau})
 \langle a_0 U_0^\ast , \ldots,  [D_\tau, \s^{-(j-1)} (a_{j-1} U_{j-1}^\ast )]_\s \cdot 
 \s^{-(j-1)} ( a_{j} U_{j}^\ast),
  \ldots,    \rangle_{D_\tau} \\ 
 &&-\,   \i_\s (\dot{D_\tau})
 \langle a_0 U_0^\ast , \ldots,  \s^{-(j-1)} (a_j U_j^\ast )\cdot [D_\tau, \s^{-j} (a_{j+1} U_{j+1}^\ast)]_\s ,
  \ldots,  \rangle_{D_\tau}\, = \\
  &&=  \sum_{k \leq j-2}  (-1)^k 
   \langle \s^4(a_0 U_0^\ast) , \ldots,  \dot{D_\tau},  \ldots , 
   [D_\tau, \s^{-(j-3)} (a_{j-1} U_{j-1}^\ast )]_\s ,  \\
   && \qquad \qquad  \qquad \qquad  \qquad \qquad  
[D_\tau^2 , \s^{-(j-1)} ( a_{j} U_{j}^\ast)]_\s, [D_\tau, \s^{-j} (a_{j+1} U_{j+1}^\ast)]_\s ,
  \ldots,  \rangle_{D_\tau}\,  \\
    &&+ \, (-1)^{j-1} 
   \langle \s^2(a_0 U_0^\ast) , \ldots , [D_\tau, \s^{-(j-3)} (a_{j-1} U_{j-1}^\ast )]_\s \cdot 
 \s^{-(j-3)} ( a_{j} U_{j}^\ast), \dot{D_\tau} , 
  \ldots,    \rangle_{D_\tau}\,  \\ 
&& -\,  (-1)^{j-1}
 \langle \s^2(a_0 U_0^\ast)  , \ldots, \s^{-(j-3)} (a_{j-1} U_{j-1}^\ast) ,\dot{D_\tau} , \\
&& \qquad \qquad \qquad \qquad \qquad \qquad
  \s^{-(j-1)} (a_j U_j^\ast )\cdot [D_\tau, \s^{-j} (a_{j+1} U_{j+1}^\ast)]_\s ,
  \ldots,  \rangle_{D_\tau} +\\
   &&+  \sum_{k \geq j}  (-1)^k
   \langle \s^4(a_0 U_0^\ast) ,\ldots , [D_\tau, \s^{-(j-5)} (a_{j-1} U_{j-1}^\ast )]_\s , 
[D_\tau^2, \s^{-(j-3)} ( a_{j} U_{j}^\ast)]_\s, \ldots \\
&&\qquad \qquad \qquad \qquad\qquad \qquad \qquad \qquad \ldots,  \dot{D_\tau}, 
  \ldots,    \rangle_{D_\tau}\, ,
   \end{eqnarray*}  
 where we have again applied Lemma \ref{forb}.  
 We focus on the last two terms, 
   \begin{eqnarray*}
  &&   \iota(\dot{D_\tau}) \langle a_0 U_0^\ast , [D_\tau, \s^{-1} (a_{1} U_{1}^\ast)]_\s ,
 \ldots,  [D_\tau, \s^{-(q-1)} (a_{q-1} U_{q-1}^\ast )]_\s \cdot 
 \s^{-(q-1)} ( a_{q} U_{q}^\ast) \rangle_{D_\tau} \\ 
 &&\qquad \qquad \qquad     -  \iota(\dot{D_\tau})
 \langle a_q U_q^\ast \cdot a_0 U_0^\ast , 
  \ldots , [D_\tau , \s^{-(q-1)} (a_{q-1} U_{q-1}^\ast)]_\s  \rangle_{D_\tau}\, \\
  && =   \langle \s^2(a_0 U_0^\ast) , \dot{D_\tau}, [D_\tau, \s^{-1} (a_{1} U_{1}^\ast)]_\s , \ldots \\
&& \qquad \qquad \qquad \qquad\qquad \qquad \qquad
\ldots,  [D_\tau, \s^{-(q-1)} (a_{q-1} U_{q-1}^\ast )]_\s \cdot 
 \s^{-(q-1)} ( a_{q} U_{q}^\ast) \rangle_{D_\tau} \\ 
 &&\qquad \qquad  - \langle \s^2(a_q U_q^\ast \cdot a_0 U_0^\ast) , \dot{D_\tau} ,
  \ldots , [D_\tau , \s^{-(q-1)} (a_{q-1} U_{q-1}^\ast)]_\s  \rangle_{D_\tau}\, + \ldots
  \end{eqnarray*}  
because one has to use
Eq.   \eqref{cyc1} to prepare them for the application of Lemma \ref{forb}, as follows:
 \begin{eqnarray*}
 && = (-1)^{(q-1)}   \langle \dot{D_\tau}, [D_\tau, \s^{-1} (a_{1} U_{1}^\ast)]_\s , 
  \ldots \\
  && \ldots , [D_\tau, \s^{-(q-1)} (a_{q-1} U_{q-1}^\ast )]_\s \cdot 
 \s^{-(q-1)} ( a_{q} U_{q}^\ast) ,  \s^{-(q-1)} ( a_0 U_0^\ast) \rangle_{\mu(U_0)D_\tau}\,  \\ 
 &&  \qquad + (-1)^q  \langle \dot{D_\tau}, [D_\tau, \s^{-1} (a_{1} U_{1}^\ast)]_\s , 
  \ldots  [D_\tau , \s^{-(q-1)} (a_{q-1} U_{q-1}^\ast)]_\s,  \\
  && \qquad \qquad  \qquad \qquad \qquad \s^{-(q-1)}(a_q U_q^\ast)
   \cdot \s^{-(q-1)}(a_0 U_0^\ast) \rangle_{\mu(U_qU_0) D_\tau}\,  .
  \end{eqnarray*}   
 In doing so, we have rescaled the operator $D_\tau$. This kind of rescaling, which
 appears every time we need to make cyclic rearrangements, prevents the
 expression \eqref{Jbco2} from vanishing.

 However, in the special case of the scaling family $\tau = t \mapsto D_t = tD$,
one can integrate from $0$ to $\infty$, replacing the ordinary integral near $0$ with its
finite part as in \cite[\S 4]{C-M1}. Also, since $D$ is invertible, the
proof of Lemma 2 in \opcit can be easily replicated to produce the necessary
estimates for the behavior of $\, J^\bullet (tD) $ and $\, \Jb^{\bullet} (tD, D)$ 
as $t \nearrow \infty$.
After integrating all the expressions involved in the
above calculation, the mismatching disappears and all cancelations that
take place in the untwisted case do occur in this case too.  Indeed, taking as an example
the first of the two terms above,  one has
\begin{eqnarray*}
  && {\rm Pf}_0 \int_\ve^{\infty}  \langle a_0 U_0^\ast, D, [D_t, \s^{-1} (a_{1} U_{1}^\ast)]_\s , \ldots \\
 &&\qquad \qquad \qquad \qquad 
 \ldots,  [D_t, \s^{-(q-1)} (a_{q-1} U_{q-1}^\ast )]_\s \cdot 
 \s^{-(q-1)} ( a_{q} U_{q}^\ast)  \rangle_{t D}\, dt \\ 
   &&=  {\rm Pf}_0 \int_\ve^{\infty}  \langle a_0 U_0^\ast, D, [D, \s^{-1} (a_{1} U_{1}^\ast)]_\s , \ldots \\
&&\qquad \qquad \qquad \qquad  \ldots , [D, \s^{-(q-1)} (a_{q-1} U_{q-1}^\ast )]_\s \cdot 
 \s^{-(q-1)} ( a_{q} U_{q}^\ast)  \rangle_{D} (t)\, \frac{dt}{t} .
    \end{eqnarray*}  
 By Eq. \eqref{cyc1} this equals
 \begin{eqnarray*}
   &&=  {\rm Pf}_0 \int_\ve^{\infty}  \langle D, [D, \s^{-1} (a_{1} U_{1}^\ast)]_\s , \ldots \\
&& \ldots  , [D, \s^{-(q-1)} (a_{q-1} U_{q-1}^\ast )]_\s \cdot 
 \s^{-(q-1)} ( a_{q} U_{q}^\ast) ,\,  \s^{-q} (a_0 U_0^\ast) \rangle_{D} (\mu(U_0)^2 t)\, \frac{dt}{t} ,
      \end{eqnarray*}  
which after the substitution $t \mapsto \mu(U_0)^{-2} t$ becomes
 \begin{eqnarray*}
    &&=  {\rm Pf}_0 \int_{ \mu(U_0)^2 \ve}^{\infty}  \langle D, [D, \s^{-1} (a_{1} U_{1}^\ast)]_\s , \ldots \\
&& \ldots  , [D, \s^{-(q-1)} (a_{q-1} U_{q-1}^\ast )]_\s \cdot 
 \s^{-(q-1)} ( a_{q} U_{q}^\ast) , \, \s^{-q} (a_0 U_0^\ast) \rangle_{D} (t)\, \frac{dt}{t} .
    \end{eqnarray*}
    \medskip
    
In this way one obtains the following {\it transgression formula}:

 \begin{lem} \label{lemtrans} For any $q \geq 0$, one has
  \begin{eqnarray}  \label{trans}
 &&{\rm Pf}_0 J^q (tD) -\, \lim_{t \nearrow \infty}J^q (tD) \, = \,\\ \notag
&&= b \left( {\rm Pf}_0 \int_\ve^{\infty} \Jb^{q-1} (tD, D) \, dt \right)\, + \,
  B \left({\rm Pf}_0 \int_\ve^{\infty} \Jb^{q+1} (tD, D) \, dt \right) .
\end{eqnarray}
\end{lem}
  
On the other hand, with virtually identical arguments as in the proof of
\cite[Proposition 2]{C-M1}, one establishes the similar vanishing result:
  \begin{eqnarray} \label{vanish}
   \lim_{t \nearrow \infty}J^q (tD) \, = \,0 \, . 
\end{eqnarray}
In particular, for $q=p$ (summability dimension), applying the $b$-boundary to
Eq. \eqref{trans}
and using the cocycle
identity \eqref{cocyid} together with the vanishing property \eqref{vanq}, one obtains
\begin{equation*}
b  \,  B  \left({\rm Pf}_0 \int_\ve^{\infty}  \Jb^{p+1} (tD, D)\, dt \right)\, = \, 
b \gimel^p (D)) \, = \, - B  \gimel^{p+2} (D)) \, = \, 0 \, .
\end{equation*}
This shows that
   \begin{eqnarray}
 \daleth^p (D) \, := \,  B \left({\rm Pf}_0 \int_\ve^{\infty}  \Jb^{p+1} (tD, D) \, dt \right) 
   \end{eqnarray}
is a cyclic cocycle. 

  \smallskip

 \begin{lem} \label{cotrans}
 The $(b, B)$-cocycle $\, \gimel^\bullet (D)\, $
 is cohomologous to the cyclic cocycle $ \daleth^p (D)$.
 \end{lem}

\proof In view of Eqs. \eqref{trans} and \eqref{vanish}, the
difference between the two cochains is a total coboundary in the
periodic cyclic complex:
  \begin{eqnarray} 
\gimel^\bullet (D) -  \daleth^p (D) \, = \,   (b+B) \left({\rm Pf}_0
 \int_\ve^{\infty} \Jb^\bullet (tD, D) \, dt  \right) .
\end{eqnarray}
\endproof
  \smallskip

  We are now ready to conclude the proof of the Ansatz for spectral triples
  twisted by scaling automorphisms.

\begin{thm} \label{cochar}
The periodic cyclic cohomology class in  $HP^\ast (\cA_G)$ of the
cocycle   $\, \gimel^\bullet (D)\, $ coincides with
 the Connes-Chern character
$Ch^\ast (\cA_G , \fH , D)$.
 \end{thm} 
 
  \proof The strategy for the proof remains the same as in \cite{C-M1, C-M2}, and 
  relies on employing the  family $D_u = D |D|^{-u} , \, u \in [0, 1]$ in order
  to construct a homotopy between the cocycle $\,  \daleth^p (D) \, $ 
  and the global cocycle $\, \tau^p_F \, $. Using the fact that each $D_u $ defines
  its own spectral triple twisted by scaling automorphisms (with character $\mu^{1-u}$), and
  with similar analytic estimates and
  algebraic manipulations as above, one establishes
  the analogue of \cite[Proposition 3]{C-M2} in the form
    \begin{eqnarray*} 
 \daleth^p(D_0) -   \daleth^p (D_1) \, = \,   
(b+B) \left( \int_0^1\,{\rm Pf}_0 \int_\ve^{\infty}  \Jb^\bullet (tD_u, D_u) \, dt  \, \, du \right)\,  .
\end{eqnarray*}
Since $D_1 = F$ and $F^2 = \id$, the cyclic cocycle
 $  \daleth^p (D_1) $ can be easily seen to coincide, up to the constant
 factor $\frac{\Gamma (\frac{p}{2}+ 1)}{2 p!}$,
  with the very cocycle $\tau^p_F$ (\cf 
 Eq. \eqref{CCh}) that
 defines the Connes-Chern character.
 \endproof

\subsection{The non-invertible case} As noted after Definition \ref{consim},
the passage from $(\cA, \fH, D)$ to the invertible double $(\wt{\cA}, \wt{\fH}, \wt{D})$
necessitates the replacement of the exact similarity condition \eqref{sim} by the
perturbed version \eqref{psim}. This does affect the twisted commutators, but only up to
higher order in the asymptotic expansion. More precisely,
\begin{equation}
[\ve^{\frac{1}{2}} \wt{D}, \td{a} \, U^\ast]_\s \, = \, 
[\ve^{\frac{1}{2}} \wt{D}, \td{a}] \,U^\ast \, + \, \ve^{\frac{1}{2}} \,a \,U^\ast \hot e_1 F_1 .
\end{equation}
At the same time,
 \begin{equation*}
\wt{D}^2 \, = \, (D^2 + \id) \hot \id \, , \qquad \text{hence} \qquad
U\, \wt{D}^2 \, U^\ast \, = \, \big(\mu (U)^2 \, D^2 \, + \, \id \big) \hot \id \, .
 \end{equation*}
Retracing the arguments leading to the expansion Eq. \eqref{asyN}, one sees that
the constant term remains unaffected by the perturbation. 
\medskip

Alternatively, one could proceed as in~\cite[\S 6.1]{Higson}, and add a compact 
operator `mass' to the Dirac Hamiltonian. Specifically, in the construction
of the invertible double, one takes
 \begin{equation*}
\wt{D} \, =\, D \hot \id + K \hot F_1\, ,  
  \end{equation*}
where   $K \in {\rm OP}^{-\infty}$ is a smoothing operator such that
\begin{equation*}
[K, D] \, = \, [D, K] \quad \text{and} \quad D^2 + K^2 \quad \text{is invertible} .
  \end{equation*}
The similarity condition is again perturbed, this time by a smoothing operator.
Since the residue integral factors through the complete symbols, the 
constant term remains of the same form as in Proposition \ref{cterm}, only with
$D$ replaced by $\wt{D}$.
\bigskip

\subsection{Application to foliations with transverse similarity structure} We conclude
by briefly indicating how one can use the above result in order to compute the index pairing
for the leaf space of a foliation with transverse similarity structure. 

A codimension $n$ foliation $\cF$ of a $N$-dimensional
manifold $V$ has a {\it transverse similarity structure} if
there exist an open cover 
$\{U_i \}_{i\in I}$ of $V$ and a family $\{h_i:U_i \ra \R^n \}_{i\in I}$  of submersions such
that $\cF\vert U_i = \{h_i^{-1}(y) ;
\,  y \in h_i(U_i)\}$ and the covering transformations 
$g_{ji} : h_i \vert U_i \cap U_j \ra h_j \vert U_j \cap U_j $ 
are given by similarities in $\Sim (n)$. Concrete 
examples of such foliations can be found in \cite{nish}, where for the case
$n = N - 1$  all nonsingular flows which admit a closed transversal (satisfying an additional
property) are in fact classified.
When $N = 3$ the notion of a transverse similarity structure to a nonsingular
flow coincides with that of a {\it complex affine structure}, treated in 
\cite{ghys} without the requirement for the existence of a closed transversal.

Given a foliation $\cF$ with a transverse similarity structure, let $\cG$ denote the smooth
\'etale groupoid associated to a complete transversal $M$ and let 
$\cA_\cG = C_c^\infty (\cG)$ (see \cite[II, \S\S8-10]{book}). The Dirac operator $D$ on $M$
defines a spectral triple twisted by similarities over the algebra $\cA_\cG$, whose
Connes-Chern character is given by the cocycle  $\, \gimel^\bullet (D) \in CC^\bullet (\cA_\cG)$.
On the other hand, let $P$ be a proper $\cG$-manifold \cite[II, \S10]{book} with compact
quotient $P / \cG$,
and let $\cD$ be a
$\cG$-invariant elliptic differential operator on $P$. By a construction
explained in \cite[\S5]{C-M} for discrete groups and in \cite[III, 7.$\gamma$]{book}
for \'etale groupoids, one associates to $\cD$ a well-defined K-theory class
$$ 
\Ind (\cD) \in K_\ast (\cA_\cG \ot \cR) \, ,
$$
where $\cR$ is the algebra of infinite matrices with rapidly decaying entries.
In the even case, this class can be represented by a difference idempotent 
$$ E_\cD - E_0 \in M_k (\cA_\cG \ot \cR) \, , \qquad k >>0 .
$$
The index pairing between the K-homology class of $D$ and the K-theory class
of $ \Ind (\cD)$ is then
computed by the pairing of their explicitly expressed Chern characters:
\begin{equation*}
<  \gimel^\bullet (D) ,  ch_\bullet (E_\cD - E_0) >  \, = \, <D, \Ind (\cD) > \,  \in \, \Z \, .
\end{equation*}


\bigskip
  
\begin{thebibliography}{01-09}

\bibitem{B-G-V} {\sc N. Berline, E. Getzler, M. Vergne},
{\bf Heat kernels and Dirac operators},
Grundlehren der Mathematischen Wissenschaften, 298. Springer-Verlag, Berlin, 1992. 
 
\bibitem{bourg} {\sc J.-P. Bourguignon, P. Goduchon}, Spineurs, op\'erateurs de Dirac
et variations,
{\it Commun. Math. Phys.} {\bf 144} (1992), 581-599.

\bibitem{ncdg} {\sc A. Connes}, Noncommutative differential geometry,
\textit{Inst. Hautes Etudes Sci. Publ. Math.} \textbf{62} (1985),
257-360.


\bibitem{book} {\sc A. Connes}, {\bf Noncommutative geometry},
 Academic Press, 1994.
 \url{ftp://ftp.alainconnes.org/book94bigpdf.pdf}
 
 \bibitem{C-M0} {\sc A. Connes, H. Moscovici},
 Transgression du caract\`ere de Chern et cohomologie cyclique,
 {\it C.R. Acad. Sci. Paris} {\bf 303} (1986), 913-918.
 
 \bibitem{C-M} {\sc A. Connes, H. Moscovici},
 Cyclic cohomology, the Novikov conjecture and hyperbolic groups,
\textit{Topology} \textbf{29} (1990), 345--388.


\bibitem{C-M1} {\sc A. Connes, H. Moscovici},
Transgression of the Chern character of finite-dimensional K-cycles,
\textit{Commun. Math.
Phys.} \textbf{155} (1993), 103-122.

\bibitem{C-M2} {\sc A. Connes, H. Moscovici},
 The local index formula in noncommutative
geometry, {\it Geom. Funct. Anal.} {\bf 5} (1995), 174-243.

\bibitem{C-M3} {\sc A. Connes, H. Moscovici}, Hopf algebras, cyclic
cohomology and the transverse index theorem, \textit{Commun. Math.
Phys.} \textbf{198} (1998), 199-246.
 
 \bibitem{C-M4} {\sc A. Connes, H. Moscovici}, Type III and spectral triples,
\textbf{Traces in Geometry,  Number Theory and Quantum Fields},
Aspects of Mathematics \textbf{E 38},
Vieweg Verlag 2008, pp. 57--71. 
\url{http://www.arxiv.org/abs/math.OA/0609703}

 \bibitem{Dave} {\sc S. Dave}, An equivariant noncommutative residue,
math.AP/061037.
\url{http://www.arxiv.org/abs/math.AP/0610371}

 \bibitem{fer} {\sc J. Ferrand},  The action of conformal
transformations on a Riemannian manifold, 
\textit{Math. Ann.}  \textbf{304} (1996), 277-291.

 \bibitem{ghys} {\sc E. Ghys},  The action of conformal
transformations on a Riemannian manifold, 
\textit{M\'em. Soc. Math. France}  \textbf{46} (1991), 123-150.

\bibitem{Var} {\sc J. M. Gracia-Bondia, J. C. Varily, H. Figueroa}, 
{\bf  Elements of Noncommutative Geometry},
{\it Birkh\"auser Advanced Texts}, Birkh\"auser Boston Inc., Boston,MA, 2001.

\bibitem{GetzSzen} {\sc E. Getzler, A. Szenes}, On the Chern character of a theta-summable
Fredholm module, \textit{J. Func. Anal.} \textbf{84} (1989), 343-357.

\bibitem{Higson} {\sc N. Higson},  The  residue index theorem of Connes and
Moscovici. {\it Surveys in noncommutative geometry},  71--126,
{\bf  Clay Math. Proc., 6}, Amer. Math. Soc., Providence, RI, 2006.
 
\bibitem{Horm} {\sc L. H\"ormander}, {\bf  The analysis of linear partial differential operators. I},
{\it Distribution theory and Fourier analysis}.  Classics in Mathematics, Springer-Verlag, Berlin, 1983.

\bibitem{JLO} {\sc A. Jaffe, A. Lesniewski, K. Osterwalder}, Quantum $K$-theory: I. The Chern
character,  \textit{Commun. Math. Phys.} \textbf{118} (1988), 1-14.


\bibitem{kulk} {\sc R. Kulkarni}, Conjugacy classes in M(n), in: 
{\bf Conformal Geometry}, R. Kulkarni and U. Pinkall (Eds.),
Aspects  Math. E 12, Vieweg, p. 41-64.

 \bibitem{nish} {\sc T. Nishimori},  A note on the classification of nonsingular flows with transverse
 similarity structures, 
\textit{Hokkaido Math. J.}  \textbf{21} (1992), 381-393.

\end{thebibliography}
\end{document}